%% file: Steering_TeX_File.tex
\documentclass{amsart}

\input{preamble}

\begin{document}

\title{Homogeneous Sasakian and 3-Sasakian Structures from the Spinorial Viewpoint}
\author{Jordan Hofmann}
\begin{abstract}
We give a spinorial construction of Sasakian and $\ts$ structures in arbitrary dimension, generalizing previously known results in dimensions $5$ and $7$. Furthermore, we obtain a complete description of the space of invariant spinors on a homogeneous $\ts$ space, and show that it is spanned by the Clifford products of invariant differential forms with a certain invariant Killing spinor. Finally, we give a basis for the space of invariant Riemannian Killing spinors on a homogeneous $\ts$ space, and determine which of these induce the homogeneous $\ts$ structure.
\end{abstract}

\maketitle


\noindent


\tableofcontents
\section{Introduction}

\input{introduction}

\textbf{Acknowledgments:} The author would like to thank Marie-Am\'{e}lie Lawn for her patience and guidance throughout the writing process, and for her feedback on initial drafts. The author would also like to thank Ilka Agricola for helpful discussions. This work was supported by the Engineering and Physical Sciences Research Council [EP/L015234/1, EP/W522429/1]; the EPSRC Centre for Doctoral Training in Geometry and Number Theory (The London School of Geometry and Number Theory); University College London; King's College London; and Imperial College London.

\input{preliminaries}
\input{main_body}

\bigskip
\bigskip

\bibliography{bibliodatabase}
\bibliographystyle{alpha}


\end{document}

%% file: preamble.tex
\usepackage{amssymb}
\usepackage{amsmath}
\usepackage{amsthm}
\usepackage{mathtools}
\usepackage{mltex}
\usepackage{parskip}
\usepackage{enumerate}
\usepackage{hyperref}
\usepackage{tcolorbox}
\usepackage{graphicx}
\graphicspath{ {images/} }
\usepackage{amsmath,microtype}
\usepackage{multicol}
\usepackage{todonotes}
\usepackage{wrapfig}
\usepackage{tikz-cd}
\usepackage[english]{babel}
\usepackage [autostyle, english = american]{csquotes}
\MakeOuterQuote{"}
\usetikzlibrary{backgrounds}
\usepackage{dynkin-diagrams}
\usepackage{booktabs}
\usepackage{blkarray}
\usetikzlibrary{babel}
\usepackage{upgreek}
\usepackage{soul}
\DeclareMathOperator{\Ad}{Ad}
\DeclareMathOperator{\ad}{ad}

\DeclareMathOperator{\Aut}{Aut}
\DeclareMathOperator{\GL}{GL}
\DeclareMathOperator{\SL}{SL}
\DeclareMathOperator{\SO}{SO}
\DeclareMathOperator{\SU}{SU}
\DeclareMathOperator{\G}{G}
\DeclareMathOperator{\Sp}{Sp}
\DeclareMathOperator{\U}{U}

\DeclareMathOperator{\Spin}{Spin}

\DeclareMathOperator{\Id}{Id}
\DeclareMathOperator{\inv}{inv}
\DeclareMathOperator{\Span}{span}

\DeclareMathOperator{\round}{round}
\DeclareMathOperator{\Alt}{Alt}

\theoremstyle{plain}
\newtheorem{theorem}{Theorem}[section]
\newtheorem{lemma}[theorem]{Lemma}
\newtheorem{proposition}[theorem]{Proposition}
\newtheorem{corollary}[theorem]{Corollary}

\theoremstyle{definition}
\newtheorem{definition}[theorem]{Definition}
\newtheorem{remark}[theorem]{Remark}

\newtheorem*{theorem*}{Theorem}

\newcommand{\Z}{\mathbb{Z}}

\newcommand{\R}{\mathbb{R}}
\newcommand{\C}{\mathbb{C}}

\newcommand{\tad}{{\text{3-}(\alpha,\delta)\text{-Sasaki}}}
\newcommand{\ts}{{\text{3-Sasakian}}}

\newcommand{\RP}{\mathbb{RP}}
\renewcommand{\:}{\colon}

\usepackage[margin=3cm]{geometry}

\newtcolorbox{mybox}
{colframe = purple!25,
  colback  = purple!10,
  coltitle = purple!20!black,  
  title    = Idea}
  
\newtcolorbox{mybox2}
{colframe = red!75,
  colback  = red!10,
  coltitle = red!20!white,  
  title    = Problem}

\newtcolorbox{mybox3}
{colframe = blue!75,
  colback  = blue!10,
  coltitle = blue!20!white,  
  title    = To be completed/added later}  
  

%% file: introduction.tex
The existence of special spinor fields on a Riemannian manifold efficiently encodes a great deal of geometric information, and is relevant, for example, to the study of immersion theory, Einstein metrics, holonomy theory, and $G$-structures, among others \cite{Friedrich_submanifold,Wang,Bar,Friedrich80,dim67}. The most extensively investigated special spinors are the real Riemannian Killing spinors, i.e. those satisfying the differential equation $\nabla^g_X\psi =\pm \frac{1}{2}X\cdot \psi$ for all vector fields $X$, whose existence places strong constraints on the geometry of the underlying manifold. Indeed, it was shown by Friedrich in \cite{Friedrich80} that Killing spinors are eigenspinors realizing the lower bound for eigenvalues of the Dirac operator, and that any manifold carrying such spinors is Einstein with scalar curvature $R=n(n-1)$. The classification of complete simply-connected manifolds with real Killing spinors was subsequently accomplished by B\"{a}r in \cite{Bar}, where it was shown that they correspond to parallel spinors (or equivalently, to a reduction of holonomy) on the metric cone. Comparing with Wang's classification of geometries carrying parallel spinors \cite{Wang}, this shows that Killing spinors are in fact somewhat rare and, beyond round spheres and isolated cases in dimensions $6$ and $7$, are carried only by (Einstein-)Sasakian and $3$-Sasakian manifolds, i.e. Riemannian manifolds $(M^{2n-1},g)$ such that the holonomy of the metric cone $(M\times \R^+, r^2 g+dr^2)$ reduces to a subgroup of $\U(n)$ or $\Sp(n/2)$; for this reason ($3$-)Sasakian manifolds play an outsized role in spin geometry. Importantly, $\ts$ manifolds $(M^{4n-1},g)$ may be considered in full generality from the perspective of spin geometry due to Kuo's result that they admit a reduction of the structure group of the tangent bundle to the simply-connected subgroup $\{1\}\times \Sp(n-1)$ of $\SO(4n-1)$ and hence are necessarily spin \cite{3Sas_structure_reduction}. Similarly, it is well-known (see e.g.\@ \cite{BG_3Sas_paper}) that simply-connected Einstein-Sasakian manifolds are spin, though Sasakian manifolds in general are not necessarily Einstein (in contrast to the well-known result of Kashiwada in the $\ts$ setting \cite{Kashiwada}). In the present article we shed new light on the correspondence between Einstein-Sasakian (and hence also $\ts$) structures and Killing spinors by giving an explicit construction of the former in terms of the latter:

\begin{theorem*}\emph{
		Let $(M,g)$ be a Riemannian spin manifold carrying a pair $\psi_1,\psi_2$ of Killing spinors (resp.\@ four Killing spinors $\psi_1,\psi_2,\psi_3,\psi_4$) for the same Killing number $\lambda \in \{\frac{1}{2},\frac{-1}{2}\}$. If the vector field $\xi_{\psi_1,\psi_2}$ defined by the equation 
		\[
		g(\xi_{\psi_1,\psi_2},X) := \Re \langle \psi_1,X\cdot \psi_2\rangle 
		\]for all $X\in TM$ has locally constant non-zero length (resp.\@ if the vector fields $\xi_{\psi_1,\psi_2}$, $\xi_{\psi_3,\psi_4}$ are orthogonal and have locally constant non-zero lengths), then, rescaling if necessary, this vector field determines a Sasakian structure on $M$ (resp.\@ these vector fields determine a $\ts$ structure on $M$). Conversely, any Einstein-Sasakian (resp.\@ $\ts$) structure on a simply-connected manifold arises by this construction.}
\end{theorem*}

Using a new argument valid in all dimensions, this theorem generalizes previous results of Friedrich and Kath in dimensions 5 and 7 \cite{FK88_french,dim5_Killing_spinors,Fried90}, which were proved by employing certain special spinorial properties occuring in these dimensions.

In the latter sections of this article we concern ourselves mainly with $\ts$ manifolds, which initially appeared over fifty years ago (see \cite{3Sas_structure_reduction,udriste69}, among others). Notable milestones in the subject include Konishi's construction of $\ts$ structures on certain principal $\SO(3)$-bundles over quaternionic K\"{a}hler manifolds of positive scalar curvature \cite{konishi75}, and the result of Boyer et al.\@ that, if the Reeb vector fields are complete, the leaf space of the induced $3$-dimensional foliation is a quaternionic K\"{a}hler orbifold \cite{BG3Sas}. Indeed, these results show that $\ts$ spaces lie between quaternionic K\"{a}hler geometries below and hyperK\"{a}hler geometries above, emphasizing the fact that they provide natural examples of interesting odd dimensional quaternionic geometries.

Previous work on $\ts$ manifolds from the spinorial perspective includes e.g.\@ \cite{Fried90,3Sasdim7}, which give a thorough and elegant accounting of the situation in dimension $7$, however until now very little is known about these spaces in higher dimensions. In Sections \ref{invarianceofstructures}, \ref{invariantspinorsandformssection}, and \ref{KSsection} of this paper we provide, in arbitrary dimension, a detailed spinorial picture of homogeneous $\ts$ spaces, which were classified by Boyer et al.\@ in \cite{BG3Sas}. Using Kuo's $\Sp(n-1)$-reduction together with the description of these spaces by Draper et al.\@ in terms of \emph{3-Sasakian data} \cite{homdata}, we apply our new invariant theoretic approach developed in \cite[Section 4]{AHLspheres} to give a classification of the invariant spinors:
\begin{theorem*}
	\emph{For a simply-connected homogeneous $3$-Sasakian manifold $(M^{4n-1}=G/H,g,\xi_i,\eta_i,\varphi_i)$, the space of invariant spinors forms an algebra under the wedge product, and is isomorphic to the algebra of invariant $\varphi_1$-anti-holomorphic differential forms:
		\[
		\Sigma_{\inv} \cong  \Lambda_{\inv}^{0,\bullet}(T_{\C}^*M),
		\]
		Furthermore, this algebra is generated by the forms $y_1:=\frac{1}{\sqrt{2}}(\xi_2+i\xi_3)$ and $\omega:=-\frac{1}{2}(\Phi_2\rvert_{\mathcal{H}}+i\Phi_3\rvert_{\mathcal{H}})$, where $\Phi_i:=g(-,\varphi_i(-))$.}
\end{theorem*}

Finally, we use this to give a complete description of the invariant Killing spinors on homogeneous $\ts$ spaces. Indeed, a partial construction of the Killing spinors on $\ts$ manifolds is given in \cite{Fried90} as sections of certain rank two subbundles of the spinor bundle, however this description can produce at most six linearly independent Killing spinors and is thus incomplete for spaces of dimension $>19$ (cf.\@ \cite{Bar}). In Section \ref{KSsection} we resolve this issue in the homogeneous setting, obtaining the following result:
\begin{theorem*}\emph{
		If $n\geq 2$ then the space of invariant Killing spinors on a simply-connected homogeneous $\ts$ manifold $(M^{4n-1}=G/H, g,\xi_i,\eta_i,\varphi_i)$ has a basis given by
		\[
		\psi_{k}:= \omega^{k+1}-i(k+1)y_1\wedge \omega^k, \qquad  -1\leq k \leq n-1,
		\]
		where we use the conventions $\omega^{-1}=0$ and $\omega^0=1$. If $n=1$ then the space of invariant Killing spinors has a basis given by $1$, $y_1$. Furthermore, if $(M,g)\ncong (S^{4n-1},g_{\round})$ then any Killing spinor is invariant.}
\end{theorem*}
As a consequence of this characterization, we deduce explicit formulas for the Killing spinors which recover the homogeneous $\ts$ strucuture via the preceding construction.

The paper is organized as follows: Section \ref{prelimssection} reviews the necessary definitions and background that will be used throughout the paper. Section \ref{constructionsection} gives a spinorial construction of Einstein-Sasakian and $\ts$ structures on simply connected manifolds, and shows that this construction recovers all such structures. Section \ref{invarianceofstructures} establishes, in the homogeneous $\ts$ setting, the compatibility of this construction with the transitive group action. Section \ref{invariantspinorsandformssection} gives a classification of the invariant spinors on homogeneous $\ts$ spaces and discusses the relationship with invariant differential forms. Section \ref{KSsection} gives explicit bases for the space of Riemannian Killing spinors carried by these spaces,
and determines which of these spinors induce the homogeneous $\ts$ structure via the construction in Section \ref{constructionsection}.

%% file: preliminaries.tex
\section{Preliminaries}\label{prelimssection}
In this section we give basic definitions and background related to the spin representation, spinors on homogeneous spaces and (homogeneous) Sasakian and $\ts$ structures. For a thorough introduction to these topics, among others, we recommend \cite{LM,BFGK,BG3Sas,BG_3Sas_paper,FriedrichBook,SasakianGeometry,homdata}.

\subsection{The Spin Representation via Exterior Forms}\label{understandingspinrep}
Throughout this paper we shall make use of the realization of the spin representation in terms of exterior forms. This realization is well-known in the context of representation theory, however its application to spin geometry and the study of spinors has not yet been widely adopted outside of \cite{Wang,AHLspheres}.  For a detailed description of this construction we refer the reader to \cite[Chap.\@ 6.1.2]{GoodmanWallach} (beware their different convention for the Clifford relation), and for its application to spinorial calculations on homogeneous spaces see \cite{AHLspheres}. We briefly recall here the basic definitions and properties insofar as they relate to this work. 

Consider $(V=\R^{2n-1},g)$ with the Euclidean inner product $g$ and the standard $g$-orthonormal basis $\{e_1,\dots, e_{2n-1}\}$. Letting $\varphi\: V\to V$ denote the almost complex structure on $(\R e_1)^{\perp}$ given by 
\[
\varphi(e_{2j})= e_{2j+1}, \quad \varphi(e_{2j+1}) = -e_{2j},
\] the complexification of $V$ can be written as a direct sum 
\begin{align}\label{VdecompLagrangians}
V^{\C} = \C_0\oplus L\oplus L', 
\end{align}
where $\C_0:=\C e_1$ and $L$ (resp. $L'$) denotes the space of $\varphi$-holomorphic (resp. $\varphi$-anti-holomorphic) vectors. Explicitly, these spaces are given by
	\begin{align}\label{LLprimespaces}
	L&:= \text{span}_{\C} \{ x_j:= \frac{1}{\sqrt{2}} (e_{2j}-ie_{2j+1})\}_{j=1}^{n-1} , \quad 
	L':= \text{span}_{\C} \{ y_j:= \frac{1}{\sqrt{2}} (e_{2j}+ie_{2j+1})\}_{j=1}^{n-1}.
\end{align}  
Letting $u_0:= ie_1$, we define an action of $V$ on the algebra $\Sigma:= \Lambda^{\bullet} L'$ of $\varphi$-anti-holomorphic forms  via
\begin{align}
	u_0:= -\Id\rvert_{\Sigma^{\text{even}}} +\Id\rvert_{\Sigma^{\text{odd}}},\quad x_j\cdot \eta &:= i\sqrt{2} \ x_j\lrcorner \eta,\quad 
	y_j\cdot \eta := i\sqrt{2} \ y_j\wedge \eta , \label{clifford_rep_formulas}
\end{align}
where $\Sigma^{\text{even}}$ and $\Sigma^{\text{odd}} $ denote the even and odd graded parts of $\Sigma=\Lambda^{\bullet} L'$. Recalling the definition of the complex Clifford algebra,
\[
\C l(V^{\C}, g^{\C}) := T(V^{\C})/(v\cdot w + w\cdot  v = -2g^{\C}(v,w)1 ),
\]
one easily verifies using the identities (5.43) in \cite{GoodmanWallach} that the action (\ref{clifford_rep_formulas}) descends to a representation of $\C l(V^{\C}, g^{\C})$ on $\Sigma$. Solving for the real orthonormal basis vectors $e_1,\dots, e_{2n-1}$ in (\ref{LLprimespaces}) gives
\begin{align*}
e_{1}= -iu_0, \quad e_{2j} = \frac{1}{\sqrt{2}} (x_j+y_j),\quad e_{2j+1} = \frac{-i}{\sqrt{2}} (y_j-x_j) \qquad \text{for all } \  j=1,\dots, n-1,
\end{align*}
and the corresponding action on $\Sigma$ (i.e. the \emph{Clifford multiplication}) is given by
\begin{align}\label{cliffordmultONB}
		e_{1}=  i\Id\rvert_{\Sigma^{\text{even}}} -i\Id\rvert_{\Sigma^{\text{odd}}},\quad e_{2j}\cdot \eta  &=  i(x_j \lrcorner \eta + y_j\wedge \eta ),\quad e_{2j+1}\cdot \eta =  (y_j\wedge\eta - x_j\lrcorner\eta ),
\end{align}
for all $\eta\in \Sigma$. It is well-known that the group $\Spin(2n-1)$ may be realized inside the Clifford algebra $\C l(V^{\C},g^{\C})$ as
\[
\Spin(2n-1) = \{ v_1 \cdot  \dotso \cdot   v_{2k} \: v_i \in \R^{2n-1}, \ ||v_i||=1 , \  k\in \mathbb{N}   \} ,
\]
and the natural extension of the Clifford multiplication to an action of $\Spin(2n-1)$ on $\Sigma$ realizes the \emph{spin representation} (see e.g.\@ \cite[Chap.\@ 1.1]{BFGK}).  
\begin{remark}
	It is possible to similarly define the Clifford multiplication and spin representation for even dimensional spaces by deleting the $\C_0$ factor in the decomposition (\ref{VdecompLagrangians}) and the corresponding operators $u_0$, $e_{1}$ in (\ref{clifford_rep_formulas}) and (\ref{cliffordmultONB}).
\end{remark}

\subsection{Spinors on Homogeneous Spaces}
Let $M=G/H$ be a reductive homogeneous space for a semisimple group $G$, and fix a reductive decomposition $\mathfrak{g}=\mathfrak{h}\oplus_{\perp} \mathfrak{m}$ which is orthogonal with respect to the Killing form on $\mathfrak{g}$. In this section we review the construction of some geometrically relevant bundles as homogeneous bundles associated to the projection $\pi \: G\to G/H$. For a more detailed introduction to reductive homogeneous spaces we refer to \cite{Arvan}, and for worked examples illustrating the process of finding invariant spinors we recommend \cite{AHLspheres} and Chapters 4.5, 5.4 of \cite{BFGK}.

First, we recall that the tangent space $T_oM$ at the origin $o:= eH$ is naturally identified with the reductive complement $\mathfrak{m}$. The other tangent spaces are then obtained by displacements of $\mathfrak{m}$ under the isometries in $G$, and there is a corresponding description of the tangent bundle as a homogeneous bundle:
\[
TM= G\times_{\Ad\rvert_H} \mathfrak{m},
\]
where $\pi\: G\to G/H$ is viewed as a principal $H$-bundle. The natural representation of $H=\text{Stab}_G(o)$ on $T_oM$ (by letting $h\in H$ act via $dh_o\: T_oM\to T_oM$) is called the \emph{isotropy representation}, and it is isomorphic to the restricted adjoint representation $\Ad\rvert_H\: H\to \GL(\mathfrak{m})$, $h\mapsto \Ad(h)\rvert_{\mathfrak{m}}$. Under this identification, an invariant Riemannian metric on $M$ corresponds to a inner product $g\:\mathfrak{m}\times \mathfrak{m}\to \R$ with the property that $\Ad\rvert_H(H) \subseteq \SO(\mathfrak{m},g) \subseteq GL(\mathfrak{m})$; for an invariant metric, the oriented frame bundle is then given as a homogeneous bundle by
\[
P_{\SO}=G\times_{\Ad\rvert_H} \SO(\mathfrak{m},g) .
\]
Suppose now that there exists a lift $\widetilde{\Ad\rvert_H}$ of the isotropy representation to the spin group, i.e.\@ a group homomorphism $\widetilde{\Ad\rvert_H}$ such that the diagram in Figure \ref{fig:hom_spin_struct} commutes.
\begin{figure}
\begin{tikzcd}
& \Spin(\mathfrak{m},g) \arrow[d, "2:1", two heads] \\
H \arrow[ru, "\widetilde{\Ad\rvert_H}", dashed] \arrow[r, "\Ad\rvert_H"'] & \SO(\mathfrak{m},g)                              
\end{tikzcd}
\caption{Homogeneous Spin Structure}
\label{fig:hom_spin_struct}
\end{figure}
Such a map induces a spin structure and spinor bundle as the homogeneous bundles
\begin{align}\label{Ginvspinstructure}
P_{\Spin}:=G \times_{\widetilde{\Ad\rvert_H}} \Spin(\mathfrak{m},g), \quad  \Sigma M:=P_{\Spin}\times_{\sigma} \Sigma = G\times_{\sigma\circ \widetilde{\Ad\rvert_H}} \Sigma, 
\end{align}
where $\sigma\: \Spin(\mathfrak{m}) \to \Aut(\Sigma)$ denotes the spin representation. Furthermore, for a connected isotropy group $H$, this was shown to be the unique $G$-invariant spin structure on $M=G/H$ by Daura Serrano, Kohn, and Lawn in \cite{invariantspinstructures} (see also the earlier works \cite{Cahen_Gutt_spin_struct_symmetric_spaces,Hirzebruch_Slodowy_elliptic_genera} for certain special cases).

In this paper we shall mainly be concerned with $\ts$ spaces (see Section \ref{prelims:Sasakian_introduction}), which are necessarily spin due to Kuo's reduction of the structure group of the tangent bundle to the (simply-connected) symplectic group of the horizontal distribution \cite[Thm.\@ 5]{3Sas_structure_reduction}. For a homogeneous $\ts$ space, denote by $\mathfrak{m}=\mathfrak{m}_{\mathcal{V}}\oplus \mathfrak{m}_{\mathcal{H}}$ the splitting into vertical and horizontal distributions. Invariance of the structure tensors implies that the image of $H$ under the isotropy representation is contained in the above reduction, i.e. 
\[
\Ad\rvert_{H} (H) \subseteq \{1\}\times  \Sp(\mathfrak{m}_{\mathcal{H}})\subseteq \SO(\mathfrak{m}),
\]
and one therefore obtains a (unique, if $H$ is connected) lifting of the isotropy representation and the associated (unique) $G$-invariant spin structure as above. Throughout the paper we will always use this spin structure when considering homogeneous $\ts$ spaces.

In light of the associated bundle construction of the spinor bundle in (\ref{Ginvspinstructure}), spinors are identified with $H$-equivariant maps $\varphi\: G\to \Sigma$, i.e. maps satisfying
\begin{align}\label{spinoridentification}
\varphi(gh) = \sigma\circ\widetilde{\Ad\rvert_H}(h^{-1})\cdot  \varphi(g) \quad \text{for all } g\in G,h\in H. 
\end{align}

\emph{$G$-invariant spinors} (which we shall refer to simply as \emph{invariant spinors} when the group $G$ is clear from context) correspond to constant $H$-equivariant maps $\varphi\: G\to \Sigma$, and we denote by $\Sigma_{\inv}\subseteq \Sigma M$ the space of such spinors. Equivalently, one sees from (\ref{spinoridentification}) that invariant spinors correspond to trivial subrepresentations of $\sigma\circ \widetilde{\Ad\rvert_H}\: H \to \GL(\Sigma)$. 

One may similarly realize the bundles of $k$-tensors and differential $k$-forms on $M$ as homogeneous bundles via
\[
\otimes^k TM = G\times_{(\Ad\rvert_H)^{\otimes k}} \mathfrak{m}^{\otimes k}, \qquad \Lambda^k TM = G\times_{\Lambda^k(\Ad\rvert_H)} \Lambda^k\mathfrak{m},
\]
and invariant sections then correspond to trivial $H$-subrepresentations of $\mathfrak{m}^{\otimes k}$ and $\Lambda^k\mathfrak{m}$ respectively. The representation-theoretic problem of finding trivial subrepresentations is approached in this paper using results from classical invariant theory, together with computer calculations in LiE (\cite{LiE}) for certain cases involving the exceptional Lie groups.

\subsection{Sasakian and $\ts$ Structures}\label{prelims:Sasakian_introduction}
Let us briefly define Sasakian and 3-Sasakian structures and discuss some of their properties, following the exposition in \cite{BG3Sas}. 
%
%
%
%
%
%
%
\begin{definition}
	A Sasakian structure on a Riemannian manifold $(M,g)$ is a unit length Killing vector field $\xi$ such that the endomorphism field $\varphi:= -\nabla^g\xi$ satisfies
	\[
	(\nabla^g_X\varphi)(Y)= g(X,Y)\xi - \eta(Y)X \quad \text{for all } X,Y\in TM  
	\]
(where $\nabla^g$ denotes the Levi-Civita connection). It is customary to call $\xi$ the \emph{Reeb vector field} and denote a Sasakian structure by $(M,g,\xi,\eta,\varphi)$, where $\eta:=\xi^{\flat}$. A Sasakian structure induces \emph{vertical} and \emph{horizontal} distributions defined by 
	\[
	\mathcal{V}:= \R \xi ,\qquad \mathcal{H}:= \ker\eta ,
	\]
	and a \emph{fundamental $2$-form} defined by
	\[
	\Phi(X,Y):= g(X,\varphi(Y)) \quad \text{for all }  X,Y\in TM.
	\]
\end{definition}
In the next proposition we collect several basic algebraic properties of Sasakian structures.
\begin{proposition}\label{sasakiantensoridentities}(Based on \cite[Prop.\@ 2.2]{BG3Sas}). If $(M,g,\xi,\eta,\varphi)$ is a Sasakian manifold, then
	\begin{align*}
	\varphi^2 &= -\Id +\eta \otimes \xi,\quad \eta(\xi)=1, \quad \varphi(\xi)=0, \quad \text{\emph{Im}}(\varphi)\subseteq \mathcal{H}, \quad d\eta = 2\Phi,\\
	0&=N_{\varphi}(X,Y) :=  [\varphi(X),\varphi(Y)] +\varphi^2[X,Y]-\varphi[\varphi(X),Y]-\varphi[X,\varphi(Y)] +d\eta(X,Y)\xi , \\
	0&= g(\varphi(X),Y) + g(X,\varphi(Y)), \quad g(\varphi(X),\varphi(Y)) = g(X,Y) - \eta(X)\eta(Y),
	\end{align*}
	for all $X,Y\in TM$.
\end{proposition}
In particular, the first identity in the preceding proposition shows that $\varphi$ is an almost complex structure on the (codimension $1$) horizontal distribution, hence a Sasakian manifold is necessarily odd-dimensional. In a similar spirit, for manifolds of dimension $4n-1$ we have the notion of a $\ts$ structure, which consists of three orthogonal Sasakian structures whose Reeb vector fields satisfy the commutator relations of the imaginary quaternions: 
\begin{definition}
	A $\ts$ structure on a Riemannian manifold $(M,g)$ consists of three Sasakian structures $(g,\xi_i,\eta_i,\varphi_i)$, $i=1,2,3$ such that the Reeb vector fields $\xi_i$, $i=1,2,3$ are orthogonal and satisfy
	\[
	[\xi_i,\xi_j]=2\xi_k
	\]
for any even permutation $(i,j,k)$ of $(1,2,3)$. It is customary to denote a $\ts$ structure by $(M,g,\xi_i,\eta_i,\varphi_i)$, omitting the "$i=1,2,3$". The \emph{vertical} and \emph{horizontal} distributions are defined by
\[
\mathcal{V}:= \text{span}_{\R}\{ \xi_i\}_{i=1}^3,\qquad \mathcal{H}:= \cap_{i=1}^3 \ker(\eta_i),
\]
and the vector fields $\xi_i$, $i=1,2,3$ are called the \emph{Reeb vector fields}. The \emph{fundamental $2$-forms} are defined by
\[
\Phi_i(X,Y):= g(X,\varphi_i(Y)) \quad   \text{for all } X,Y\in TM.
\]
\end{definition}
In addition to the identities in Proposition \ref{sasakiantensoridentities}, the tensors defining a $\ts$ structure satisfy certain "pseudo-quaternionic" compatibility relations:
\begin{proposition} (Based on \cite[Eqn.\@ (2.4)]{BG3Sas} and \cite[Eqn.\@ (1.5)]{3str}). If $(M,g,\xi_i,\eta_i,\varphi_i)$ is a $\ts$ manifold, then
\begin{align*}
\varphi_i&= \varphi_j\circ \varphi_k - \eta_k\otimes \xi_j = -\varphi_k\circ \varphi_j + \eta_j\otimes \xi_k \\
\varphi_i(\xi_j)&=-\varphi_j(\xi_i)= \xi_k ,  \qquad \eta_i=\eta_j\circ \varphi_k = -\eta_k\circ\varphi_j ,
\end{align*}
for any even permutation $(i,j,k)$ of $(1,2,3)$.
\end{proposition}
In particular, one sees from the preceding proposition that $\varphi_i$, $i=1,2,3$ determine a quaternionic structure on the (codimension $3$) horizontal distribution, hence a $\ts$ manifold necessarily has dimension congruent to $3$ modulo $4$. This quaternionic structure on $\mathcal{H}$ naturally leads to a certain type of local frame which is well-adapted to the geometry:
\begin{definition}
	Let $(M,g,\xi_i,\eta_i,\varphi_i)$ be a $\ts$ manifold. A local frame $e_1,\dots, e_{4n-1}$ of $TM$ is called \emph{adapted} if 
	\[
	e_i=\xi_i, \qquad e_{4p+i}=\varphi_i(e_{4p}) , \qquad \text{ for all } i=1,2,3, \ \  p = 1,\dots, n-1. 
	\]
\end{definition}

%
%
%
%
%
%
%
%
%
%
%
%
%
%
%
%
Finally, for $\ts$ manifolds we recall that there is a particularly useful choice of metric connection adapted to the geometry: the so-called \emph{canonical connection} (see \cite[Section 4]{3str}, where this connection is introduced for the more general class of $\tad$ manifolds). We review its main properties in the following proposition: 
\begin{proposition}(Based on \cite[Section 4]{3str}).\label{canonicalconnectiondefinition}
	For a $\ts$ manifold $(M,g,\varphi_i,\xi_i,\eta_i)$, the canonical connection $\nabla$ is the unique metric connection with skew torsion such that
\begin{align*}
\nabla_X\varphi_i &= -2 ( \eta_k(X) \varphi_j -\eta_j(X) \varphi_k) \qquad \text{for all } X\in TM.
\end{align*}
The derivatives of the other structure tensors are
\[
\quad \nabla_X\xi_i = -2 ( \eta_k(X) \xi_j -\eta_j(X) \xi_k) , \quad 
\nabla_X\eta_i = -2 ( \eta_k(X) \eta_j -\eta_j(X) \eta_k),
\]
and the torsion 3-form is given by $T= \sum_{i=1}^3 \eta_i\wedge d\eta_i$.
\end{proposition}

\subsection{Homogeneous $\ts$ Spaces}
First, we recall Boyer, Galicki, and Mann's classification of homogeneous $\ts$ spaces:
\begin{theorem}\label{hom3Sas_classification}
	(Based on \cite[Thm.\@ C]{BG3Sas}). The homogeneous $\ts$ spaces $(M^{4n-1}=G/H,g)$ are precisely
	\begin{align*}
		& S^{4n-1} \cong \frac{\Sp(n)}{\Sp(n-1)}, \quad  \mathbb{RP}^{4n-1}\cong \frac{\Sp(n)}{\Sp(n-1)\times \mathbb{Z}_2}, \quad \frac{\SU(n+1)}{\text{\emph{S}}( \U(n-1) \times \U(1)) }, \\
		& \frac{\SO(n+3)}{\SO(n-1)\times \Sp(1)}, \quad \frac{\G_2}{\Sp(1)}, \quad \frac{\text{\emph{F}}_4}{\Sp(3)}, \quad  \frac{\text{\emph{E}}_6}{\SU(6)}, \quad \frac{\text{\emph{E}}_7}{\Spin(12)}, \quad \frac{\text{\emph{E}}_8}{\text{\emph{E}}_7},
	\end{align*}
	where the permissible values of $n$ are as follows:
	\begin{table}[h!] 
		\centering
		\begin{tabular}{ |c||c|c|c| }
			\hline
			$G$	&  $\Sp(n)$ & $\SU(n+1)$ & $\SO(n+3)$ \\
			\hline
			n &   $n\geq 1$ & $n\geq 2 $ & $n\geq 4 $  \\
			\hline
		\end{tabular}
	\end{table}
\end{theorem}
This classification was obtained by proving that any homogeneous $\ts$ space fibers over a Wolf space with a finite list of possibilities for the fiber, and then using the classification of Wolf spaces in \cite{wolf_spaces}. Recently, Goertsches et al.\@ obtained a new proof of the classification using root systems of complex simple Lie algebras to construct homogeneous $\ts$ spaces \cite{Leander_roots}. Previously, and also from the algebraic point of view, the invariant connections on homogeneous $\ts$ spaces were studied in detail by Draper et al.\@ in \cite{homdata}. Importantly, they gave a characterization of these spaces in terms of purely Lie theoretic data called \emph{3-Sasakian data}, which we recall in Theorem \ref{DOP_3Sas_data}. Before this, we briefly recall in Proposition \ref{Nom_maps_prop} the notion of \emph{Nomizu maps}, which are the algebraic objects encoding invariant connections on homogeneous vector bundles. For our purposes it will only be necessary to consider Nomizu maps in the context of metric connections on the tangent bundle and their lifts to spinor bundle, however it should be noted that the theory is valid in greater generality. For a modern exposition of the topic we refer the reader to the forthcoming book \cite{ANT_book_principal_fibre_bundles}, and the classical references are the original papers of Nomizu \cite{Nomizumap} and Wang \cite{Wangconnections} and the introductory books of Kobayashi and Nomizu \cite{KN1,KN2}.
\begin{proposition}\label{Nom_maps_prop}
	Let $(M=G/H,g)$ be a Riemannian homogeneous space with origin $o:=eH$ and reductive decomposition $\mathfrak{g}=\mathfrak{h}\oplus \mathfrak{m}$, and let $\mathfrak{so}(\mathfrak{m})$ denote the special orthogonal Lie algebra with respect to the $\Ad(H)$-invariant inner product on $\mathfrak{m}$ induced by $g$ under the standard identification $\mathfrak{m}\cong T_oM$. There is a one-to-one correspondence between metric connections $\nabla$ on $TM$ and linear maps $\Uplambda\: \mathfrak{m}\to \mathfrak{so}(\mathfrak{m})$ intertwining the actions of $\Ad(H)$. If $\omega$ is an invariant tensor or differential form, then its covariant derivative is given by
	\[
	(\nabla_{\widehat{X}}\omega)_o  = \Uplambda(X)\omega_o \quad \text{for all } X \in \mathfrak{m}, 
	\]
	where $\widehat{X}:= \frac{d}{dt}\rvert_{t=0} L_{\exp_G(tX)}$ is the fundamental vector field associated to $X\in \mathfrak{m}$ and $\Uplambda(X)\in \mathfrak{so}(\mathfrak{m})$ acts on $\omega_o$ by the natural extension of the action on $\mathfrak{m}$.
\end{proposition}
\begin{theorem}\label{DOP_3Sas_data}(Based on \cite[Def.\@ 4.1, Thm.\@ 4.3]{homdata}). Let $M^{4n-1}=G/H$ be a homogeneous space with connected isotropy group $H$, satisfying the following properties:
	\begin{enumerate}[(i)]
		\item $\mathfrak{g}$ is compact;
		\item $\mathfrak{g}$ is simple and there is a $\Z_2$-graded decomposition $\mathfrak{g}=\mathfrak{g}_0\oplus \mathfrak{g}_1$ such that $\mathfrak{g}_0=\mathfrak{sp}(1)\oplus \mathfrak{h}$\emph{;}
		\item There exists an $\mathfrak{h}^{\C}$-module $U$ of complex dimension $2(n-1)$ such that $\mathfrak{g}_1^{\C}\cong \C^2 \otimes U$ as a module for $\mathfrak{g}_0^{\C}\cong \mathfrak{sp}(1)^{\C}\oplus \mathfrak{h}^{\C}$, where $\C^2$ is the standard representation of $\mathfrak{sp}(1)^{\C} \cong \mathfrak{sl}(2,\C)$. 
	\end{enumerate}
	Then there is a homogeneous $\ts$ structure $(\xi_i,\eta_i,\varphi_i)$ on $M=G/H$ determined by the tensors
	\begin{align*}
		\xi_1&:=  \begin{pmatrix}
			i&0\\0&-i
		\end{pmatrix}, \quad \xi_2:= \begin{pmatrix}
			0&-1\\1&0
		\end{pmatrix} , \quad \xi_3:= \begin{pmatrix}
			0&-i\\-i&0
		\end{pmatrix},\quad \varphi_i =\frac{1}{2}\ad(\xi_i)\rvert_{\mathfrak{sp}(1)} + \ad(\xi_i)\rvert_{\mathfrak{g}_1},\\
		g&:= -\frac{1}{4(n+1)}\kappa\rvert_{\mathfrak{sp}(1)\times\mathfrak{sp}(1)} -\frac{1}{8(n+1)}\kappa\rvert_{\mathfrak{g}_1\times \mathfrak{g}_1}, 
	\end{align*}
	where $\kappa$ denotes the Killing form of $\mathfrak{g}$. Furthermore, the Nomizu map of the Levi-Civita connection $\nabla^g$ is given by
	\begin{align}\label{LCNomizumap}
		\Uplambda^g(X)Y = \begin{cases}
			\frac{1}{2}[X,Y]_{\mathfrak{m}} & \text{ if $X,Y\in\mathfrak{sp}(1)$ or $X,Y\in \mathfrak{g}_1$},\\
			0 & \text{ if $X\in\mathfrak{sp}(1)$, $Y\in \mathfrak{g}_1$},\\
			[X,Y]_{\mathfrak{m}} & \text{ if $X\in\mathfrak{g}_1$, $Y\in\mathfrak{sp}(1)$},
		\end{cases}
	\end{align}
	where subscript $\mathfrak{m}$ denotes projection onto the reductive complement $\mathfrak{m}:=\mathfrak{sp}(1)\oplus \mathfrak{g}_1$.
\end{theorem}
Indeed, they proved that all simply-connected homogeneous $\ts$ spaces can be constructed from $\ts$ data, and they gave an explicit description of the data in each case which will be extremely useful for our purposes.

Finally, we recall the result of Agricola, Dileo, and Stecker that, in the homogeneous case, the Nomizu map of the canonical connection (see Proposition \ref{canonicalconnectiondefinition}) takes a simple form:
\begin{proposition}(Based on \cite[Prop.\@ 4.2.1]{hom3alphadelta}). For a homogeneous $\ts$ space, the Nomizu map $\Uplambda$ of the canonical connection $\nabla$ is given by
	\begin{align}\label{canonicalNomizumap} 
		\Uplambda(X)Y = \begin{cases}
			-[X,Y] & \text{ if $X\in \mathcal{V}$,} \\
			\ 0& \text{ if $X\in\mathcal{H}$},
		\end{cases}
	\end{align}
	for all $X,Y\in\mathfrak{m}$.
\end{proposition}

%% file: main_body.tex
\section{The ($3$-)Sasakian Structures Induced by Killing Spinors}\label{constructionsection}
In this section we expain how to construct Sasakian and $\ts$ structures from Riemannian Killing spinors, generalizing the constructions of Friedrich and Kath in dimensions $5$ and $7$ \cite{FK88_french,dim5_Killing_spinors,Fried90}. By considering the rank two subbundles $E_i^-\subseteq \Sigma M$ defined in \cite{Fried90}, we also show that all Sasakian and $\ts$ structures on connected, simply-connected Einstein-Sasakian or $\ts$ manifolds arise from this construction.
\begin{definition}
	A spinor $\psi\in\Gamma(\Sigma M)$ is called a (Riemannian) Killing spinor if it satisfies
	\[
	\nabla_X^g \psi = \lambda X\cdot \psi
	\]
	for all $X\in TM$. We shall refer interchangeably to Riemannian Killing spinors and Killing spinors, and, unless otherwise stated, we will only consider real Killing spinors for the Killing numbers $\lambda=\pm \frac{1}{2}$ (such a normalization can always be achieved by an appropriate rescaling of the metric: compare \cite[Cor.\@ 1 on p.\@ 21]{BFGK} with \cite[Thm.\@ 7 on p.\@ 29]{BFGK}).
\end{definition}
Generalizing to arbitrary dimension the $1$-form and dual vector field considered in \cite{FK88_french} (also those considered in \cite[Section 5]{Fried90}, \cite[Chapter 4.4]{BFGK}, and, by letting the two spinors in Definition \ref{associated_tensors_definition} be multiples of each other, the $1$-form considered in \cite[Section 4]{dim5_Killing_spinors}), we make the following definition:
\begin{definition}\label{associated_tensors_definition}
Let $(M,g)$ be a spin manifold. Given a pair of spinors $\psi_1,\psi_2$, the \emph{associated $1$-form} $\eta_{\psi_1,\psi_2}$ and its metric dual $\xi_{\psi_1,\psi_2}$ (the \emph{associated vector field}) are defined by
\begin{align}\label{inducedvectorfield}
\eta_{\psi_1,\psi_2}(X) := \Re \langle \psi_1,  \ X\cdot \psi_2\rangle, \qquad \xi_{\psi_1,\psi_2}:= \eta_{\psi_1,\psi_2}^{\sharp}, \qquad \text{for all } X\in TM,
\end{align}
where $\langle\cdot,\cdot\rangle $ denotes the usual Hermitian metric on the spinor bundle and $\Re$ is the real part. 
%
%
We also define the \emph{associated endomorphism field} $\varphi_{\psi_1,\psi_2}$ by
\begin{align}\label{inducedendomorphism}
\varphi_{\psi_1,\psi_2}(X)  := -\frac{1}{2} (X\lrcorner d\eta_{\psi_1,\psi_2})^{\sharp}, \qquad \text{for all } X\in TM.
\end{align}
\end{definition}
\begin{remark}
	Observe that the operator $\Re$ in (\ref{inducedvectorfield}) is unnecessary when $\psi_1\in \R i \psi_2$, since skew-symmetry of the Clifford multiplication with respect to the Hermitian product $\langle \ , \ \rangle$ ensures that $\langle i\psi,X\cdot \psi\rangle$ is purely real. Thus, by setting $\psi_1=i\psi_2$ in (\ref{inducedvectorfield}), our definition recovers the $1$-form and vector field considered in \cite[Section 4]{dim5_Killing_spinors} and \cite[Chapter 1.5]{BFGK}. We also note that it is possible to choose non-vanishing spinors $\psi_1,\psi_2$ such that the associated vector field $\xi_{\psi_1,\psi_2}$ is identically zero. This is necessarily the case when $\psi_1=\psi_2$, for example.
\end{remark}
As we will frequently encounter the tensors (\ref{inducedvectorfield}) and (\ref{inducedendomorphism}) and their derivatives, we summarize here some relevant well-known facts:
\begin{lemma}\label{realpartcommutes} Let $(M,g)$ be a Riemannian spin manifold with spinor bundle $\Sigma M$, and denote by $\langle \ ,\ \rangle$ the usual Hermitian scalar product on the fibers of $\Sigma M$.
	\begin{enumerate}[(i)]
	\item Differentiation of the Hermitian product $\langle \ , \ \rangle$ commutes with $\Re$, i.e.
	\[
	X(\Re\langle \varphi,\psi\rangle ) = \Re(X\langle \varphi,\psi\rangle)
	\]
	for all vector fields $X$ and spinors $\varphi,\psi$.
	\item For any spinor $\psi \in\Gamma(\Sigma M)$ and vector fields $X,Y\in TM$ we have 
	\[
	\Re \langle X \cdot \psi, Y\cdot \psi \rangle = g(X,Y)\ ||\psi||^2.
	\]
	\item If $\psi$ is a Riemannian Killing spinor then the length function $||\psi||$ is constant on each connected component of $M$.
	\end{enumerate}
\end{lemma} 
\begin{proof}
	These facts are easy to check and, in some cases, appear frequently throughout the literature (see e.g.\@ \cite{BFGK,GKSEinstein,GKSspheres}); in the interest of comprehensiveness we briefly recall the proofs here.
	\begin{enumerate}[(i)]
		\item One easily calculates:
	\begin{align*}
		X(\Re\langle \varphi,\psi\rangle) &= \frac{1}{2}X(\langle \varphi,\psi\rangle+ \overline{\langle \varphi,\psi\rangle})   =  \frac{1}{2}X(\langle \varphi,\psi\rangle + \langle \psi,\varphi\rangle   ) \\
		&= \frac{1}{2}\left( \langle \nabla_X^g\varphi,\psi\rangle + \langle \varphi, \nabla_X^g\psi\rangle + \langle \nabla_X^g\psi,\varphi\rangle + \langle \psi, \nabla_X^g\varphi\rangle \right) \\
		&= \frac{1}{2}\left(    \langle \nabla_X^g\varphi,\psi\rangle + \overline{\langle \nabla_X^g\varphi,\psi\rangle} +  \langle \varphi, \nabla_X^g\psi\rangle + \overline{\langle \varphi, \nabla_X^g\psi\rangle }   \right)  \\
		&= \Re \langle \nabla_X^g\varphi,\psi\rangle + \Re \langle \varphi,\nabla_X^g\psi\rangle = \Re(X\langle \varphi, \psi\rangle).
	\end{align*}
	\item Using the skew-symmetry of the Hermitian product with respect to Clifford multiplication, we calculate:
	\begin{align*}
		\Re \langle X\cdot \psi,Y\cdot \psi\rangle &= -\Re \langle \psi, X\cdot Y\cdot \psi\rangle =\Re  \langle \psi, Y\cdot X\cdot \psi\rangle  +\Re \langle \psi, 2g(X,Y)\psi \rangle   \\
		&= \Re \overline{\langle Y\cdot X\cdot \psi, \psi\rangle} + 2g(X,Y)\ ||\psi||^2 = \Re \langle Y\cdot X\cdot \psi, \psi\rangle +2g(X,Y)\ ||\psi||^2\\
		&= -\Re \langle X\cdot \psi,Y\cdot \psi\rangle + 2g(X,Y)\ ||\psi||^2 ,
	\end{align*}
and the result follows.
\item Let $\psi$ be a Riemannian Killing spinor for the Killing number $\lambda$. Differentiating the square of the norm with respect to any $X\in TM$ and using skew-symmetry of the Hermitian product with respect to Clifford multiplication gives:
\begin{align*}
	X||\psi||^2 &=  X\langle \psi, \psi \rangle =  \langle \nabla_X^g \psi,\psi\rangle +  \langle \psi, \nabla_X^g\psi \rangle = \lambda \langle X\cdot \psi,\psi \rangle + \lambda\langle \psi,X\cdot \psi\rangle = 0.
\end{align*}
	\end{enumerate}
\end{proof}

Since one of the defining conditions of a Sasakian structure involves the exterior derivative of the Reeb $1$-form, we also need an identity relating $d\eta_{\psi_1,\psi_2}$ to the spinors $\psi_1,\psi_2$. Generalizing the identity calculated by Friedrich and Kath in the proof of \cite[Thm.\@ 2]{dim5_Killing_spinors}, we have:
\begin{lemma}\label{detalemma}
	If $(M,g)$ is a Riemannian spin manifold carring a pair of Killing spinors $\psi_1,\psi_2$ (not necessarily linearly independent) for the same Killing number $\lambda$, then the exterior derivative of the $1$-form $\eta_{\psi_1,\psi_2}$ is given by
	\begin{align}\label{deta_squaring_construction_formula}
		d\eta_{\psi_1,\psi_2}(X,Y) = 2\lambda \Re \langle \psi_1, (Y\cdot X - X\cdot Y)\cdot \psi_2\rangle  \quad \text{for all } X,Y\in TM.
	\end{align}
\end{lemma}
%
%
%
%
%
%
%
%
\begin{remark}\label{squaring_construction_sasakian_tensors_remark}
	We recall here the well-known \emph{squaring construction}, by which a pair of spinors $\psi_1,\psi_2$ (not necessarily linearly independent) gives rise to differential forms $\omega_{(k)}$ of every degree $k\geq 0$:
	\begin{align} \label{squaringconstructionformula}
	\omega_{(k)}(X_1,\dots, X_k) := \Re \langle (X_1\wedge \dots \wedge X_k)\cdot \psi_1, \psi_2\rangle  	\quad \text{for all } X_1,\dots, X_k \in TM.
	\end{align}
	For our purposes it is important to note that $\eta_{\psi_1,\psi_2}= -\omega_{(1)}$ by definition, and that $d\eta_{\psi_1,\psi_2}= 4\lambda \omega_{(2)}$ by the preceding lemma.   
\end{remark}
We now describe the relationship between the length of the vector field $\xi_{\psi_1,\psi_2}$ and the kernel of $d\eta_{\psi_1,\psi_2}$:
\begin{lemma}\label{Chap_Hof22:length_function_lemma}
	If $\psi_1, \psi_2$ are Killing spinors (not necessarily linearly independent) on $(M,g)$ for the same Killing number $\lambda$, then the length function $\ell:=||\xi_{\psi_1,\psi_2}||$ is locally constant if and only if $\xi_{\psi_1,\psi_2}\lrcorner d\eta_{\psi_1,\psi_2}=0$.
\end{lemma}
\begin{proof} One easily sees from the calculation on \cite[p.\@ 30]{BFGK} that $\xi_{\psi_1,\psi_2}$ is a Killing vector field (actually they prove this only for the case $\psi_1=i\psi_2$, but the same argument works in general by using Lemma \ref{realpartcommutes}(i)). Therefore we have
	\begin{align}\label{Chap_Hof22:Killing_vector_identity}
		\frac{1}{2}d\eta_{\psi_1,\psi_2}(X,Y) = g(\nabla_X^g \xi_{\psi_1,\psi_2},Y) \quad \text{ for all } X,Y\in TM
	\end{align}
	(see e.g.\@ \cite[p.\@ 64]{blair_contact_manifolds_book}, noting the slightly different conventions), and the result then follows from the calculation	
	\begin{align*}
		X(\ell^2)&=  Xg(\xi_{\psi_1,\psi_2},\xi_{\psi_1,\psi_2}) = 2g(\nabla_X \xi_{\psi_1,\psi_2},\xi_{\psi_1,\psi_2}) = d\eta_{\psi_1,\psi_2}(X,\xi_{\psi_1,\psi_2}). 
	\end{align*}
\end{proof}
For a Sasakian structure $(g,\xi,\eta,\varphi)$ one always has $\xi\lrcorner d\eta =0$, so we see from the preceding lemma that it is necessary for $\ell$ to be locally constant (and non-zero) in order for the tensors $\frac{1}{\ell}\xi_{\psi_1,\psi_2}, \frac{1}{\ell}\eta_{\psi_1,\psi_2},\frac{1}{\ell}\varphi_{\psi_1,\psi_2}$ to constitute a Sasakian structure. We now arrive at the first main result of this section, which shows that for Killing spinors this condition is also sufficient:
\begin{theorem}\label{spinortoACM}
	Suppose that $(M,g)$ is a Riemannian spin manifold carrying a pair $\psi_1,\psi_2$ of Killing spinors (not necessarily linearly independent) for the same Killing number $\lambda \in \{\frac{1}{2},\frac{-1}{2}\}$, and suppose furthermore that $\ell:=||\xi_{\psi_1,\psi_2}||$ is locally constant and non-zero. Then, the tensors 
	\[
	\xi:=\frac{1}{\ell}\xi_{\psi_1,\psi_2},\quad \eta:=\frac{1}{\ell} \eta_{\psi_1,\psi_2},\quad \varphi:= \frac{1}{\ell} \varphi_{\psi_1,\psi_2}
	\] determine a Sasakian structure on $(M,g)$. 
\end{theorem}
\begin{proof}
	It is clear that $||\xi||=1$ by construction, and we have already noted in the proof of Lemma \ref{Chap_Hof22:length_function_lemma} that $\xi_{\psi_1,\psi_2}$ is a Killing vector field (as $\ell$ is constant on each connected component of $M$ we then have that $\xi=\frac{1}{\ell}\xi_{\psi_1,\psi_2}$ is a Killing vector field as well). Furthermore, it is clear from the definitions (\ref{inducedvectorfield}) and (\ref{inducedendomorphism}) that if $\xi$ determines a Sasakian structure on $(M,g)$ then $\eta$ and $\varphi$ are the associated $1$-form and endomorphism field. The assertion that $\xi$ determines a Sasakian structure has a very elegant proof in terms of conformal Killing forms, i.e.\@ differential forms $\Theta \in \Lambda^k T^*M$ satisfying
	\begin{align}\label{Chap_Hof22:conformal_killing_form_eqn}
		\nabla^g_X \Theta &= \frac{1}{k+1} X\lrcorner d\Theta - \frac{1}{n-k+1} X^*\wedge \delta\Theta \quad \text{for all } X \in TM,
	\end{align}
	where $\delta\:\Lambda^k T^*M \to \Lambda^{k-1} T^*M$ denotes the codifferential and $X^*$ is the $1$-form dual to $X$ via the Riemannian metric\footnote{The author would like to thank Prof.\@ Uwe Semmelmann for pointing out this argument.}. Indeed, from Remark \ref{squaring_construction_sasakian_tensors_remark} we have $d\eta=\frac{1}{\ell}d\eta_{\psi_1,\psi_2}=  \frac{4\lambda}{\ell} \omega_{(2)}$, thus by \cite[Prop.\@ 2.2]{Sem03_conformal_killing_forms} $d\eta$ is a conformal Killing form (using the fact that $\psi_1,\psi_2$ are Killing spinors, hence twistor spinors--see e.g.\@ \cite{BFGK}). The existence of the Killing spinors $\psi_1,\psi_2$ for the Killing number $\lambda\in \{\frac{1}{2},\frac{-1}{2}\}$ implies that $(M^n,g)$ is Einstein with scalar curvature $4n(n-1)\lambda^2=n(n-1)$ \cite{Friedrich80}, and the result then follows by \cite[Prop.\@ 3.4]{Sem03_conformal_killing_forms}.
\end{proof}
In fact, we shall prove in the next theorem that, \emph{any} Einstein-Sasakian structure on a simply-connected manifold arises from this construction. To that end, we recall the bundles
	\begin{align}\label{Eplusminusspaces}
	E_{\pm} := \{ \psi \in \Sigma M\: (\pm 2\varphi(X) +\xi\cdot X - X\cdot \xi) \cdot \psi =0 \ \   \text{for all } X\in TM\}
\end{align}
considered by Friedrich and Kath in \cite{Fried90}. They showed that, on a simply-connected Einstein-Sasakian spin manifold, these bundles admit a global basis consisting of Riemannian Killing spinors for the constants $\mp \frac{1}{2}$ (see the proof of \cite[Thm.\@ 1]{Fried90}). For our purposes, it is also important to recall that $\text{rank}( E_-)  \geq 1$ (see the proof of Theorem 1 in \cite[Chapter 4.2]{BFGK}, noting that the roles of $E_{+}, E_-$ are reversed here due to our convention for the Clifford multiplication). Using these bundles, we prove:
\begin{theorem}\label{Sasakian_converse}
	If $(M^{2n-1} ,g,\xi,\eta,\varphi)$ is a simply-connected Einstein-Sasakian manifold, then the Sasakian structure $(\xi,\eta,\varphi)$ arises from the preceding construction.
\end{theorem}
\begin{proof}
	Let $\psi\in\Gamma(E_-)$ be a (non-trivial) Killing spinor for the Killing number $\lambda = \frac{1}{2}$, and assume without loss of generality that $||\psi||=1$ (if not, by Lemma \ref{realpartcommutes}(iii), we may rescale it to unit length using locally constant functions). Defining $\psi':=-\xi\cdot \psi$, we note that $\psi'\in\Gamma(E_-)$, since $\xi$ anti-commutes in the Clifford algebra with the operators $(-2\varphi(X)+\xi\cdot X-X\cdot \xi)$ defining the bundle $E_-$. Furthermore, $\psi'$ is a Killing spinor due to the calculation
	\begin{align*}
		\nabla^g_X\psi' &= -\nabla_X^g(\xi\cdot \psi) = -(\nabla^g_X\xi)\cdot \psi - \xi\cdot \nabla_X^g\psi = \varphi(X) \cdot \psi -\frac{1}{2}\xi \cdot X\cdot \psi  \\
		&=  \frac{1}{2}(\xi\cdot X - X\cdot \xi)\cdot \psi	-\frac{1}{2}\xi\cdot X\cdot \psi = -\frac{1}{2}X\cdot \xi\cdot \psi = \frac{1}{2}X\cdot \psi',
	\end{align*}
	where we have used the identity $\varphi= -\nabla^g\xi$ and the condition defining $E_-$. Using Lemma \ref{realpartcommutes}(ii), we have
	\[
	\eta_{\psi,\psi'}(X) = \Re \langle \psi,X\cdot \psi'\rangle = -\Re \langle \psi,X\cdot \xi\cdot \psi\rangle = g(X,\xi) \ ||\psi||^2 = \eta(X)
	\]
	for any $X\in TM$, hence $\eta_{\psi,\psi'}=\eta$ (and $\xi_{\psi,\psi'}=\xi$, $\varphi_{\psi,\psi'} = \varphi$).
	The proof for the case $\psi \in \Gamma(E_+)$, $\lambda=-\frac{1}{2}$ is similar. 
\end{proof}
\begin{remark}
	Viewing Sasakian structures as the odd-dimensional version of K\"{a}hler structures, the bundles $E_{\pm}$ and the preceding construction are readily seen to parallel certain behaviour observed by Kirchberg in the K\"{a}hler setting \cite{KirchbergKlaus-Dieter1986Aeft}. Indeed, for a K\"{a}hler spin manifold $(M^{2m},g,J,\omega )$, Kirchberg proved that the spinor bundle splits under Clifford multiplication by the K\"{a}hler $2$-form $\omega= g(\cdot, J(\cdot))$ into eigenbundles 
	\[ \Sigma M = \Sigma_0 \oplus \dots\oplus  \Sigma_m, \qquad \text{rank}_{\C} \Sigma_r = \binom{m}{r},  \]    
	where $\omega$ acts on $\Sigma_r$ by $i(m-2r)\Id$. In particular, \cite[Prop.\@ 1]{KirchbergKlaus-Dieter1986Aeft} implies that if $\psi\in \Gamma(\Sigma_0)$ then $J(X)\cdot \psi =  i X\cdot \psi$ for all $X\in TM$. On the other hand, injectivity of Clifford multiplication by real tangent vectors implies that this relation uniquely defines $J$, so the data of the K\"{a}hler structure is encoded by the rank-$1$ bundle $\Sigma_0$ (and there is a similar relationship between $\Sigma_m$ and $-J$). Alternatively, it follows from Lemma \ref{realpartcommutes}(ii) that if $\psi \in \Gamma(\Sigma_0)$ then the squaring construction applied to the spinors $\psi_1:=-i\psi$, $\psi_2 := \psi$ gives $\omega_{(2)} = ||\psi||^2 \omega$, and hence also recovers the K\"{a}hler structure (analogously for $\Sigma_m$). This closely mirrors the picture for Einstein-Sasakian manifolds $(M,g,\xi,\eta,\varphi)$, where sections of the bundle $E_+$ (resp.\@ $E_-$) are easily seen to satisfy $\varphi(X)\cdot \psi = X\cdot \xi\cdot \psi$ (resp.\@ $\varphi(X)\cdot \psi = - X\cdot \xi \cdot \psi$) for all horizontal vector fields $X$. It follows from these relations that the squaring construction applied to $\psi_1:= -\xi\cdot \psi $, $\psi_2:= \psi$ gives $\omega_{(2)} = ||\psi||^2 \Phi $ (resp.\@ $\omega_{(2)} = -||\psi||^2\Phi$), and the Reeb vector field is similarly obtained via $\omega_{(1)} = ||\psi||^2 \xi^{\flat}$, recovering the Einstein-Sasaki structure. Indeed, this shows that certain non-zero constant length sections of $E_{\pm}$ are sufficient to recover the data of the Einstein-Sasakian structure, however in practice the easiest choice is non-zero Killing spinors as in the proof of Theorem \ref{Sasakian_converse}. For a more representation-theoretic perspective on recovering the Einstein-Sasakian structure from certain choices of Killing spinors, see e.g.\@ \cite[Section 3.4]{physics_KS_paper}, noting that the forms $\omega_{(k)}$ obtained using the squaring construction are sometimes referred to in the literature as \emph{spinor bilinears}.     
\end{remark}
As corollaries of the preceding two theorems, we also obtain analogous construction and uniqueness results for $\ts$ manifolds:
\begin{theorem}\label{spinorto3structure}
	Suppose $(M,g)$ is a Riemannian spin manifold carrying Killing spinors $\psi_1,\psi_2,\psi_3,\psi_4$ (not necessarily linearly independent) for the same Killing number $\lambda\in \{\frac{1}{2},\frac{-1}{2}\}$. If $\xi_{\psi_1,\psi_2}$ and $\xi_{\psi_3,\psi_4}$ are orthogonal vector fields with locally constant non-zero length, then the two Sasakian structures induced by Theorem \ref{spinortoACM} determine a $\ts$ structure on $(M,g)$. 
\end{theorem}
\begin{proof}
	This follows from Theorem \ref{spinortoACM} and the fact that two Sasakian structures with orthogonal Reeb vector fields uniquely determine a $\ts$ structure (see e.g.\@ \cite[p.\@ 556]{Fried90} for the construction of the structure tensors of the third Sasakian structure in terms of the other two).
\end{proof}
To prove a uniqueness result analogous to Theorem \ref{Sasakian_converse} for $\ts$ manifolds, we consider Friedrich and Kath's natural generalization
	\begin{align}\label{Eplusminusspaces_3Sas}
E_i^-:= \{ \psi \in \Sigma M\: (- 2\varphi_i(X) +\xi_i\cdot X - X\cdot \xi_i) \cdot \psi =0  \ \  \text{for all } X\in TM\}, \quad i=1,2,3
\end{align}
of the bundles (\ref{Eplusminusspaces}) to the $\ts$ setting (see the proof of \cite[Thm.\@ 6]{Fried90}). Indeed, performing the argument from the proof of Theorem \ref{Sasakian_converse} for each of the three Sasakian structures individually yields:
%
%
%
%
%
%
%
%
%
\begin{theorem}\label{3Sasakian_converse}
	If $(M^{4n-1}, g,\xi_i,\eta_i,\varphi_i)$ is a simply-connected $\ts$ manifold, then the $\ts$ structure arises from the preceding construction.
\end{theorem}
Let us briefly compare our results with those of Friedrich and Kath in the $5$- and $7$-dimensional setting. These can be found in \cite{FK88_french}, \cite[Section 4]{dim5_Killing_spinors}, \cite[Sections 5, 6]{Fried90}, and also appeared subsequently in Chapters 4.3 and 4.4 of the book \cite{BFGK}. In dimension $5$, their construction uses the fact that the spin representation of $\Spin(5)$ is transitive on the unit sphere in the spinor module (which does not hold in dimension $>9$ \cite{MontgomerySamelson43}). Given a non-zero spinor $\psi $, this allows them to arrange a particular choice of frame in which a unique unit length solution $\xi\in TM^5$ to the equation $\xi\cdot \psi = i\psi$ is readily apparent. For a non-zero Killing spinor $\psi$, and under an appropriate normalization of the scalar curvature, they then show that this vector field $\xi$ determines a Sasakian structure. Similarly, given an orthonormal pair of real Killing spinors $\psi_1,\psi_2\in \Gamma(\Sigma_{\R}M^7)$ in dimension $7$, Friedrich and Kath use the orthogonal decomposition $\Sigma_{\R} M^7= \R \psi_1 \oplus_{\perp} (TM^7\cdot \psi_1)$ of the real spinor bundle to find a unique unit length vector field $\xi$ satisfying $\xi\cdot \psi_1=\psi_2$, which goes on to become the Reeb vector field of a Sasakian structure. This decomposition of the spinor bundle occurs as a coincidence in dimension $7$, and fails in higher dimensions since the dimension of the spinor module grows much faster than the dimension of the manifold. In both cases, Friedrich and Kath note that their vector field $\xi$ is dual to the $1$-form defined essentially by (\ref{inducedvectorfield}), but there is no mention of the fact that this $1$-form can be taken as the starting point to perform such a construction in arbitrary dimension, as we have done here. 

Aside from the preceding comments about dimension, we also note that our results are slightly different in spirit: Friedrich and Kath prove that a Killing spinor in dimension $5$ (resp.\@ an orthonormal pair of Killing spinors in dimension $7$) defines a specific unit length vector field which is in fact the Reeb vector field of a Sasakian structure. On the other hand, our Theorem \ref{spinortoACM} makes no assumption about linear independence or orthogonality of the Killing spinors, but requires the induced vector field $\xi_{\psi_1,\psi_2}$ to have locally constant positive length; this is then shown in Theorem \ref{Sasakian_converse} to be a reasonable assumption in the sense that simply-connected Einstein-Sasakian manifolds always carry Killing spinors $\psi_1,\psi_2$ such that $ ||\xi_{\psi_1,\psi_2}||=1$, so no cases are lost by imposing this. %
%
The similarities and differences are much the same when comparing \cite[Section 6]{Fried90} to our Theorem \ref{spinorto3structure}.

\section{Invariance of Spinors and their Associated ($3$-)Sasakian Structures}\label{invarianceofstructures}
Given the relationship described above it is natural to ask whether, on a homogeneous manifold, invariance of an Einstein-Sasakian or $\ts$ structure implies invariance of the associated Killing spinors and vice versa. One already sees from \cite[Remark 4.43]{AHLspheres} that the homogeneous sphere $S^{4n-1}= \frac{\Sp(n)\U(1)}{\Sp(n-1)\U(1)}$ equipped with the round metric (which is invariant and carries an invariant Einstein-Sasakian structure) admits non-invariant Killing spinors $\psi\in\Gamma(E_1^-)$. However, it turns out that if one applies Theorem \ref{spinortoACM} to a pair of invariant Killing spinors then the resulting Sasakian structure must also be invariant, as we prove in this section. This suggests that an invariant spinor is a more fundamental geometric object than an invariant ($3$-)Sasakian structure, capturing more of the homogeneity data of the space. To begin, we collect several well-known and straightforward Clifford algebra commutator identities:
\begin{lemma}
	For any $X \in \R^k$ and $\theta \in\Lambda^p \R^k$ $(p\geq 0)$, the identity
	\begin{align}\label{cliffalgidentityformvector}
		\theta \cdot X - X\cdot \theta = ((-1)^p+1) \ X\lrcorner \theta + ((-1)^p -1) \ X\wedge \theta
	\end{align}
	holds in the Clifford algebra.
\end{lemma}
\begin{proof}
	The Clifford algebra identities
	\begin{align*}
		X\cdot \theta = (X\wedge \theta) - (X\lrcorner\theta),\qquad \quad  \theta\cdot X= (-1)^p[X\wedge \theta +X\lrcorner\theta]	
	\end{align*}
	appear as Equations (1.4) in \cite[Chapter 1.2]{BFGK}, and the result follows by subtracting them.
\end{proof}
Considering a $2$-form $T=\sum_{i<j} T_{ij}e_i\wedge e_j \in \Lambda^2\R^k \cong \mathfrak{so}(k)$ and its spin lift $\widetilde{T}=\frac{1}{2}\sum_{i<j}T_{ij}e_i\cdot e_j$, one immediately sees:
\begin{corollary}
	Let $T\in \mathfrak{so}(k)$ be a skew-symmetric linear transformation and $\widetilde{T}\in\mathfrak{spin}(k)$ its spin lift under the Lie algebra isomorphism $\mathfrak{so}(k)\cong \mathfrak{spin}(k)$. Then, for any $X\in \R^k$, the identity
	\begin{align}\label{cliffalgidentity}
		\widetilde{T}\cdot X  - X\cdot \widetilde{T} = T(X)
	\end{align}
	holds in the Clifford algebra. 
\end{corollary}
This formula also easily generalizes for the commutator of a 2-form with a form of arbitrary degree:
\begin{lemma}
	Let $T\in \mathfrak{so}(k)$ be a skew-symmetric linear transformation and $\widetilde{T}\in\mathfrak{spin}(k)$ its spin lift under the Lie algebra isomorphism $\mathfrak{so}(k)\cong \mathfrak{spin}(k)$. Then, for any $\theta \in \Lambda^p \R^k $ $(p\geq 0)$, the identity
	\begin{align}\label{cliffalgidentityforms}
		\widetilde{T} \cdot \theta - \theta \cdot \widetilde{T} = T(\theta),
	\end{align}
	holds in the Clifford algebra, where $T(\theta)$ refers to the standard action of $\mathfrak{so}(k)$ on $\Lambda^p \R^k $. 
\end{lemma}
\begin{proof}
	It suffices to prove the result for $T=e_i\wedge e_j$ and $\theta = e_{l_1}\wedge\dots \wedge e_{l_p}$. We calculate:{\small
		\begin{align*}
			\widetilde{T}\cdot \theta  - \theta\cdot \widetilde{T} &= \frac{1}{2} \left(   e_i\cdot e_j\cdot e_{l_1}\cdot \dotso \cdot e_{l_p} -  e_{l_1}\cdot \dotso \cdot e_{l_p} \cdot e_i\cdot e_j \right) \\
			&= \begin{cases}   0  & i,j\in \{ l_1,\dots, l_p\} \text{ or }  i,j  \notin  \{l_1,\dots, l_p\}, \\
				e_i\cdot e_j \cdot \theta  & \text{otherwise}    \end{cases} \\
			&= T(\theta).
	\end{align*}}
\end{proof}
Next, we show that the Clifford product of an invariant vector or differential form with an invariant spinor is again invariant:
\begin{lemma}\label{invformtimesinvspinor}
	If $(M=G/H,g)$ is a homogeneous spin manifold carrying an invariant spinor $\psi$, then $\theta\cdot \psi \in \Sigma $ is invariant for any invariant form $\theta\in \Lambda_{\inv}^p \mathfrak{m}$ $(p\geq 0)$.
\end{lemma}
\begin{proof}
	For any isotropy operator $h\in\ad(\mathfrak{h})\rvert_{\mathfrak{m}}\subseteq  \mathfrak{so}(\mathfrak{m})$, it follows from (\ref{cliffalgidentityforms}) that
	\begin{align*}
		\widetilde{h} \cdot \theta\cdot \psi - \theta \cdot \widetilde{h}\cdot \psi = h(\theta)\cdot \psi.
	\end{align*} 
	Invariance of $\psi$ and $\theta$ gives $\widetilde{h}\cdot \psi=0 = h(\theta)$, hence $	\widetilde{h} \cdot \theta\cdot \psi = 0$ as desired.
\end{proof}
With the preceding lemmas, it is easy to prove that invariant Killing spinors induce invariant ($3$-)Sasakian structures via the construction from Section \ref{constructionsection}: 
\begin{theorem}
	If $(M=G/H, g)$ is a Riemannian homogeneous spin manifold carrying a pair $\psi_1,\psi_2$ of invariant spinors, then the associated tensors $\xi_{\psi_1,\psi_2}$, $\eta_{\psi_1,\psi_2}$, and $\varphi_{\psi_1,\psi_2}$ are also invariant. In particular, if $\psi_1,\psi_2$ are invariant Killing spinors for the same Killing number $\lambda\in\{\frac{1}{2},\frac{-1}{2} \}$, and $\xi_{\psi_1,\psi_2}$ has locally constant positive length $\ell>0$, then the induced Sasakian structure $(\frac{1}{\ell} \xi_{\psi_1,\psi_2}, \frac{1}{\ell}\eta_{\psi_1,\psi_2}, \frac{1}{\ell} \varphi_{\psi_1,\psi_2})$ is invariant.
\end{theorem}
\begin{proof}
	We show that the tensors $\xi_{\psi_1,\psi_2}, \eta_{\psi_1,\psi_2}, \varphi_{\psi_1,\psi_2}$ are invariant. Using (\ref{cliffalgidentity}) and invariance of $g$, $\langle \ ,\ \rangle$, $\psi_1$, and $\psi_2$, we calculate:
	\begin{align*}
		g([h,\xi_{\psi_1,\psi_2}],X) &= -g(\xi_{\psi_1,\psi_2},[h,X])   =  -\Re \langle \psi_1,[h,X]\cdot \psi_2\rangle \\
		&=  -\Re \langle \psi_1,\big(\widetilde{\ad(h)\rvert_{\mathfrak{m}}}\cdot X - X\cdot \widetilde{\ad(h)\rvert_{\mathfrak{m}}}\big)\cdot \psi_2\rangle = - \Re \langle \psi_1, \widetilde{\ad(h)\rvert_{\mathfrak{m}}}\cdot X\cdot \psi_2\rangle \\
		&= \Re\langle \widetilde{\ad(h)\rvert_{\mathfrak{m}}}\cdot \psi_1,X\cdot \psi_2\rangle =0
	\end{align*}
	for all $h\in\mathfrak{h}$, $X\in \mathfrak{m}$, hence $\xi_{\psi_1,\psi_2}$ and $\eta_{\psi_1,\psi_2}=\xi_{\psi_1,\psi_2}^{\flat}$ are invariant. Invariance of $\varphi_{\psi_1,\psi_2}$ then follows from (\ref{inducedendomorphism}), completing the proof.
\end{proof}
By the same argument, one also obtains the analogous result in the $\ts$ setting:
\begin{theorem}
	If $(M=G/H, g)$ is a Riemannian homogeneous spin manifold carrying invariant Killing spinors $\psi_1,\psi_2,\psi_3,\psi_4$ for the same Killing number $\lambda\in\{\frac{1}{2},\frac{-1}{2}\}$, and such that $\xi_{\psi_1,\psi_2}$ and $\xi_{\psi_3,\psi_4}$ are orthogonal and have locally constant positive length, then the induced $\ts$ structure is invariant.
\end{theorem}
\section{Invariant Differential Forms and Spinors}\label{invariantspinorsandformssection}
Expanding upon the work \cite{homdata}, in this section we describe the invariant $\varphi_1$-(anti-)holomorphic differential forms on homogeneous $\ts$ spaces. We also describe the invariant spinors carried by these spaces, and the relationship between the forms and spinors. We would like to emphasize that this approach exploits the exterior form viewpoint of the spin representation, which greatly simplifies calculations and allows one to easily prove spinorial results for spaces of arbitrary dimension.

\begin{remark}
	Before discussing the invariant forms and spinors, we comment briefly about the homogeneous $\ts$ space $\RP^{4n-1}\cong \frac{\Sp(n)}{\Sp(n-1)\times \Z_2} $. The isotropy group $\Sp(n-1)\times \Z_2$ is not connected, leading to non-uniqueness of lifts of the isotropy representation (see Figure \ref{fig:hom_spin_struct}), and, consequently, non-uniqueness of homogeneous spin structures.
	With this in mind, we exclude $\RP^{4n-1}$ from consideration for the remainder of the article by requiring that our homogeneous $\ts$ spaces be simply-connected; it is the only non simply-connected space in Theorem \ref{hom3Sas_classification}. More details about the algebraic structure of this space may be found in \cite[Rmk.\@ 4.8]{homdata}, and for a detailed discussion of the invariant spinors on homogeneous projective spaces we refer the reader to the forthcoming article \cite{AHprojective}.

\end{remark}

%
%
%
%
%
In order to prove the first major result of this section, we will make use of the First Fundamental Theorems of Invariant Theory for the classical complex simple Lie groups, which can be found e.g.\@ in \cite{FultonHarris1991,tensorFFTs}; we will use the formulations presented in \cite{tensorFFTs} as these are more suited to our purposes. We will also need the description of the exterior powers of the standard representation of $\SO(n,\C)$ as highest weight modules (see e.g.\@ \cite[Chap.\@ 5.5.2]{GoodmanWallach}). We summarize these results in the following three theorems:
\begin{theorem}\label{fft_theorem}
	(Based on the First Fundamental Theorems in \cite[Section 5]{tensorFFTs}). Let $\SO(n,\C)$ act on $\C^n$ by its standard representation and, if $n=2l$ is even, let $\Sp(2l,\C)$ also act by its standard representation. Denote by $e_1,\dots, e_n$ (resp.\@ $e_1^*,\dots, e_n^*$) the standard basis for $\C^n$ (resp.\@ the images of the standard basis vectors under the isomorphism $\C^n \simeq (\C^n)^*$ given by the non-degenerate bilinear form defining the group), and let 
	\[
	\mathcal{T}:= 	T(\C^{n} \oplus (\C^{n})^*)
	\] 
	denote the algebra of tensors on $\C^{n} \oplus (\C^{n})^*$, with the natural algebra multiplication given by concatenation of tensors. The subalgebra of invariant tensors for the two groups are described up to mutations (i.e.\@ permutations of the tensor factors) as follows:
	\begin{enumerate}[(i)]
		\item \textul{FFT for $\SO(n,\C)$:} The space $\mathcal{T}^{\SO(n,\C)}$ of invariant tensors is the $\C$-span of all mutations of tensor products of flips of the tensors
		\[
		\det:= \sum_{\sigma\in S_n} \text{\emph{sign}}(\sigma) e_{\sigma(1)}\otimes \dots \otimes e_{\sigma(n)}, \qquad I:= \sum_{i=1}^n e_i\otimes e_i^*,
		\]
		where a flip means applying the isomorphism $\C^{n}\simeq (\C^{n})^*$, $e_i\mapsto e_i^*$ or its inverse to one of the tensor factors.
		\item \textul{FFT for $\Sp(2l,\C)$:} The space $\mathcal{T}^{\Sp(2l,\C)}$ of invariant tensors is the $\C$-span of all mutations of tensor products of $p$, $p^*$, $I$, where
		\[
		p:= \sum_{i=1}^l (e_i\otimes e_{l+i} -e_{l+i}\otimes e_i ), \qquad I:= \sum_{i=1}^n e_i\otimes e_i^*.
		\]
	\end{enumerate}
\end{theorem}

\begin{remark}
	Note that our $\mathcal{T}$ is defined slightly differently than the tensor algebra considered in \cite{tensorFFTs}; we don't require all the covariant and contravariant factors to be collected, i.e.\@ we don't require elements of $\mathcal{T}$ to lie in $T(\C^n) \otimes T((\C^n)^*)$. To compensate for this, our notion of mutations includes \emph{all} permutations of the tensor factors, not just those that permute covariant factors and contravariant factors amongst themselves (as was the definition in \cite{tensorFFTs}). It is easy to see that our mutations intertwine the group action, and that invariant elements of $\mathcal{T}$ are generated by invariant elements of $T(\C^n) \otimes T((\C^n)^*)$ by taking linear combinations of mutations.   
\end{remark}

\begin{theorem}\label{GWinvarianttheoremSLn}
	(Based on \cite[Thm.\@ 5.5.11]{GoodmanWallach}). Denote by $\omega_1,\dots, \omega_{n-1}$ the fundamental weights of $\SL(n,\C)$, and $\Lambda^r \C^n$ the $r$\textsuperscript{th} exterior power of the standard representation. The representation $\Lambda^r\C^n$ is irreducible for all $r=1,\dots,n$, with highest weight $\omega_r$ for $r=1,\dots,n-1$.
\end{theorem}
\begin{theorem}\label{GWinvarianttheorem}
	(Based on \cite[Thm.\@ 5.5.13]{GoodmanWallach}). For $n=2l$ or $2l+1$, let $\omega_1,\dots \omega_l$ denote the fundamental weights of $\SO(n,\C)$, and $\Lambda^r \C^n$ the $r$\textsuperscript{th} exterior power of the standard representation.
	\begin{enumerate}[(i)]
		\item For $n=2l+1\geq 3$: The representation $\Lambda^r\C^n$ is irreducible for all $r=1,\dots,n$, with highest weight $\omega_r$ for $r=1,\dots,l-1$ and $2\omega_l$ for $r=l$.
		\item For $n=2l\geq 4$: The representation $\Lambda^r\C^n$ is irreducible for $r=1,\dots, l-1$, with highest weight $\omega_r$ for $r=1,\dots,l-2$ and $\omega_{l-1}+\omega_l$ for $r=l-1$. The representation $\Lambda^l\C^n$ splits as the direct sum of two irreducible representations with highest weights $2\omega_{l-1}$ and $2\omega_l$.
	\end{enumerate}
\end{theorem}

We begin with an easy result regarding the algebra of $\Sp(n)$-invariant horizontal complex tensors for the classical case of the $\ts$ round sphere:

\begin{proposition}\label{invtensors}
	The algebra $T_{\inv}(\mathfrak{m}_{\mathcal{H}}^{\C}) $ of invariant horizontal complex tensors on the $\ts$ sphere $(S^{4n-1}=\frac{\Sp(n)}{\Sp(n-1)},g_{\round})$ is generated, up to mutations, by its degree $2$ elements.
\end{proposition}
\begin{proof}
	The isotropy algebra is $\mathfrak{h}=\mathfrak{sp}(n-1)$, and complexifying the isotropy representation gives $\mathfrak{m}^{\C} \simeq 3\C \oplus 2\C^{2n-2}$, where $\C^{2n-2}$ denotes the standard representation of $\mathfrak{h}^{\C}= \mathfrak{sp}(2n-2,\C)$. Using the fact that $\C^{2n-2}$ is self-dual as an $\mathfrak{h}^{\C}$-representation, the space of horizontal complex tensors is
	\begin{align}
		T(\mathfrak{m}_{\mathcal{H}}^{\C})  \simeq  T(2\C^{2n-2}) \simeq  T(\C^{2n-2} \oplus (\C^{2n-2})^*) , \label{Chap_Hof22:round_sphere_tensor_calculation}
	\end{align}
	and the result follows from Theorem \ref{fft_theorem}(ii).

\end{proof}
As a consequence, we are able to deduce a description of the $\Sp(n)$-invariant forms on the round sphere:
\begin{corollary}\label{formscor}
	For any $k\geq 0$, the space $\Lambda^k_{\inv}\mathfrak{m}$ of invariant $k$-forms on the $\ts$ sphere $(S^{4n-1}=\frac{\Sp(n)}{\Sp(n-1)},g_{\round})$ is spanned by wedge products of invariant $1$- and $2$- forms. Explicitly, the algebra $\Lambda_{\inv}^{\bullet}\mathfrak{m}$ of invariant forms is spanned by elements of the form 
	\[
	\tau_{\epsilon_1,\epsilon_2,\epsilon_3,a_1,a_2,a_3} :=  \eta_1^{\epsilon_1}\wedge \eta_2^{\epsilon_2}\wedge \eta_3^{\epsilon_3} \wedge (\Phi_1\rvert_{\mathcal{H}})^{a_1} \wedge (\Phi_2\rvert_{\mathcal{H}})^{a_2} \wedge (\Phi_3\rvert_{\mathcal{H}})^{a_3} , 
	\] 
	where $\epsilon_1,\epsilon_2,\epsilon_3 \in \{0,1\}$, $a_1,a_2,a_3\in \Z_{\geq 0}$.
\end{corollary}
\begin{proof}
	Any degree $2$ horizontal tensor decomposes uniquely into symmetric and skew-symmetric parts according to
	\[
	\mathfrak{m}_{\mathcal{H}} \otimes \mathfrak{m}_{\mathcal{H}} \simeq S^2(\mathfrak{m}_{\mathcal{H}}) \oplus \Lambda^2(\mathfrak{m}_{\mathcal{H}}),
	\]
	and this decomposition holds as representations of $H=\Sp(n-1)$. By Proposition \ref{invtensors}, we then have that $T_{\inv}(\mathfrak{m}_{\mathcal{H}}^{\C})$ is generated, up to mutations, by its symmetric and skew-symmetric degree $2$ elements:
	\[
	T_{\inv} (\mathfrak{m}^{\C}_{\mathcal{H}})/\text{mut.} \simeq  \bigoplus_{k\geq 0} ((\mathfrak{m}^{\C}_{\mathcal{H}}\otimes \mathfrak{m}^{\C}_{\mathcal{H}})_{\inv})^{\otimes k}  \simeq \bigoplus_{k\geq 0} \left( S^2_{\inv}(\mathfrak{m}^{\C}_{\mathcal{H}}) \oplus \Lambda^2_{\inv}(\mathfrak{m}^{\C}_{\mathcal{H}})\right)^{\otimes k} .
	\]
	But any tensor containing a factor from $S^2_{\inv}(\mathfrak{m}^{\C}_{\mathcal{H}})$, or any mutation of such a tensor, is symmetric in at least two positions and therefore lies in the kernel of the projection onto the skew-symmetric tensors. In particular, this shows that
	\[
	\Lambda^{\bullet}_{\inv} (\mathfrak{m}^{\C}_{\mathcal{H}}) \simeq  \bigoplus_{k\geq 0}  \Lambda^k\left(\Lambda^2_{\inv}(\mathfrak{m}^{\C}_{\mathcal{H}})\right)
	\]
	(see Case 3(i) in the proof of Proposition \ref{horizholomorphicforms} below for more details about the required properties of the skew-symmetrization map), so the algebra of complex horizontal differential forms is generated by its degree $2$ elements. For $U=\C^{2n-2}$ the standard representation of $\mathfrak{h}^{\C} = \mathfrak{sp}(2n-2,\C)$ we have $\dim_{\C} \Lambda^2_{\inv}U = 1$ (see e.g.\@ \cite[p.\@ 834]{homdata}), hence 
	\[
	\dim_{\C} \Lambda^2_{\inv}(\mathfrak{m}^{\C}_{\mathcal{H}}) = \dim_{\C} \Lambda^2_{\inv}(U\oplus U ) = 2\dim_{\C} \Lambda^2_{\inv}U + \dim_{\C} (U\otimes U)^{\mathfrak{h}^{\C}} = 2\cdot 1 + 1 = 3
	\]
	(since $U$ is a self-dual $\mathfrak{h}^{\C}$-representation), and it follows that $\Lambda^2_{\inv}(\mathfrak{m}^{\C}_{\mathcal{H}})$ is spanned by $\Phi_1\rvert_{\mathcal{H}},\Phi_2\rvert_{\mathcal{H}}$, and $\Phi_3\rvert_{\mathcal{H}}$ (viewed as complex forms). In particular, these generate the complex algebra $\Lambda^{\bullet}_{\inv}(\mathfrak{m}^{\C}_{\mathcal{H}})$ and also the real subalgebra $\Lambda^{\bullet}_{\inv}(\mathfrak{m}_{\mathcal{H}})$. The result then follows by noting that the isotropy representation acts trivially in the vertical directions. 
\end{proof}
In fact, we expect that Corollary \ref{formscor} should hold for all homogeneous $\ts$ spaces (with an additional generator $\Phi_0\in \Lambda_{\inv}^2 (\mathfrak{m}_{\mathcal{H}})$ for the case $G=\SU(n+1)$), however the arguments for the remaining cases become much more difficult (due to reducibility of $\mathfrak{m}_{\mathcal{H}}$ in the $\SU(n+1)$ case, the appearance, in certain degrees, of extra candidates for invariant tensors in the $G=\SO(n+3)$ case, and the lack of invariant-theoretic tools for the exceptional cases). It seems likely that other methods would be needed to prove such a result in general. For this reason we now prove a somewhat weaker result, but one which can be shown in all cases and which will nonetheless be sufficient for the purpose of finding the spaces of invariant spinors:
\begin{proposition}\label{horizholomorphicforms}
	If $(M=G/H, g, \xi_i,\eta_i,\varphi_i)$ is a simply-connected homogeneous $\ts$ space, then the algebras $\Lambda_{\inv}^{\bullet,0 }(\mathfrak{m}_{\mathcal{H}}^{\C}) $, $\Lambda_{\inv}^{0, \bullet }(\mathfrak{m}_{\mathcal{H}}^{\C}) $ of invariant horizontal $\varphi_1$-(anti-)holomorphic forms are generated by their degree $2$ elements.
\end{proposition}

\begin{proof}
	Employing the basis for $\mathfrak{sp}(1)=\Span_{\R}\{\xi_1,\xi_2,\xi_3 \}$ given in \cite[Eqn.\@ (17)]{homdata}, the almost complex structure $\varphi_1\rvert_{\mathcal{H}} = \ad(\xi_1) $ acts on $(\mathfrak{m}_{\mathcal{H}})^{\C}=(\mathfrak{g}_1)^{\C}$ by $i\Id $ on $(1,0) \otimes U$ and by $-i\Id $ on $(0,1)\otimes U$, hence the $\varphi_1$-holomorphic (resp.\@ $\varphi_1$-anti-holomorphic) horizontal cotangent bundles are given by $\Lambda^{1,0}(\mathfrak{m}_{\mathcal{H}}^{\C}) = (1,0)\otimes U \simeq U$ (resp.\@ $\Lambda^{0,1}(\mathfrak{m}_{\mathcal{H}}^{\C}) = (0,1)\otimes U \simeq U$). Therefore we have $\Lambda^k U \simeq \Lambda^{k,0}(\mathfrak{m}_{\mathcal{H}}^{\C})\simeq \Lambda^{0,k}(\mathfrak{m}_{\mathcal{H}}^{\C})$ for all $k\geq 0$, so it suffices to consider the invariant exterior forms on $U$. We use this approach to treat all the cases for $G$ individually:
	
	\textul{Case 1: $G=\Sp(n)$.} Here the complexified isotropy group is $\Sp(n-1)^{\C} \cong \Sp(2n-2,\C)$ and we have $U=\C^{2n-2}$ (the standard representation). It then follows from Theorem \ref{fft_theorem}(ii) that $\Lambda^k_{\inv} U \simeq \Lambda^k_{\inv} \C^{2n-2} $ is generated by alternating tensor powers (i.e.\@ wedge powers) of the $2$-form $p\in \Lambda^2\C^{2n-2}$ stabilized by $\Sp(2n-2,\C)$.

	\textul{Case 2: $G=\SU(n+1)$.} The isotropy group in this case is $H=S(\U(n-1)\times \U(1))$, and we consider separately the cases $n>2$ and $n=2$. When $n>2$ we have $U= \C^{n-1}\oplus (\C^{n-1})^*$, with the action of $\mathfrak{h}^{\C} \cong \mathfrak{sl}(n-1,\C) \oplus \mathfrak{u}(1)^{\C}$ on $U$ by the standard (resp.\@ dual of the standard) representation of $\mathfrak{sl}(n,\C)$ on $\C^{n-1}$ (resp.\@ $(\C^{n-1})^*$), and by the action of $1\in \mathfrak{u}(1)^{\C} \cong \C  $ given by 
	\begin{align}\label{u1action_n_geq_3_case}
		1\cdot v = \left(1+\frac{n-1}{2}\right)v ,\qquad  1\cdot v' =  -\left(1+\frac{n-1}{2}\right)v'
	\end{align}
	for all $v\in \C^{n-1}$, $v'\in (\C^{n-1})^*$ (see \cite[Section 4.5]{homdata}). We then have
	\begin{align}\label{Case2eqn1}
		\Lambda^k U &\simeq  \Lambda^k(\C^{n-1}\oplus (\C^{n-1})^*) \simeq  \bigoplus_{p+q=k} (\Lambda^p\C^{n-1} ) \otimes (\Lambda^q (\C^{n-1})^*),
	\end{align}
	and examining the action of $\mathfrak{u}(1)^{\C}$ in (\ref{u1action_n_geq_3_case}) shows that an element in one of the summands on the right hand side of (\ref{Case2eqn1}) is $\mathfrak{u}(1)^{\C}$-invariant if and only if it has the same number of $\C^{n-1}$ and $(\C^{n-1})^*$ factors. By Theorem \ref{GWinvarianttheoremSLn}, the $\SL(n-1,\C)$-modules $\Lambda^p\C^{n-1}$ and $\Lambda^q\C^{n-1}$ are irreducible and non-isomorphic unless $p=q$, and it then follows from (\ref{Case2eqn1}) that
	\begin{align}\label{ExteriorUdimensions}
		\dim_{\C} \Lambda^k_{\inv} (U) = \begin{cases*} 1 & \text{if $k$ is even and $k\leq \dim_{\C} U$,} \\ 0 & \text{otherwise.} \end{cases*}
	\end{align}
	In particular one checks in a basis that $\omega_{1,0}:= (\Phi_2\rvert_{\mathcal{H}}-i\Phi_3\rvert_{\mathcal{H}})$ (resp.\@ $\omega_{0,1}:=(\Phi_2\rvert_{\mathcal{H}}+i\Phi_3\rvert_{\mathcal{H}})$) is an element of $\Lambda_{\inv}^{2,0 }(\mathfrak{m}^{\C}_{\mathcal{H}})$ (resp.\@ $\Lambda_{\inv}^{0,2}(\mathfrak{m}^{\C}_{\mathcal{H}})$), and that the top power of $\omega_{1,0}$ (resp.\@ $\omega_{0,1}$) is a $\varphi_1$-holomorphic (resp.\@ $\varphi_1$-anti-holomorphic) volume form. It follows that lower powers are non-zero, hence they span the $1$-dimensional spaces of invariant $\varphi_1$-(anti-)holomorphic forms in the relevant dimensions. The argument for $n=2$ is similar, except one only needs to consider the action of $\mathfrak{h}^{\C} \cong \mathfrak{u}(1)^{\C}$ via (\ref{u1action_n_geq_3_case}).

	\textul{Case 3: $G=\SO(n+3)$.}
	The isotropy group in this case is $H= \SO(n-1)\times \Sp(1)$, and it is shown in \cite[Section 4.3]{homdata} that $U=\C^{n-1}\times \C^{n-1}$ with the following action of $H^{\C}$: identify each element $(a,b) \in U$ with $[a\rvert b]\in\text{Mat}_{(n-1)\times 2}(\C)$ (i.e.\@ the matrix with first column $a\in\C^{n-1}$ and second column $b\in \C^{n-1}$), and let $\SO(n-1,\C)$ act by left multiplication on the columns of $[a\rvert b]$ and $\Sp(1)^{\C}\cong \Sp(2,\C)$ by right multiplication on the rows. This gives decompositions $U\simeq \C^{n-1}\oplus \C^{n-1}$ as $\SO(n-1,\C)$-modules and $U\simeq (\C^2)^{\oplus (n-1)}$ as $\Sp(2,\C)$-modules (direct sums of the standard representation in both cases). Letting $e_1,\dots e_{n-1}$ (resp.\@ $e_1^*,\dots,e_{n-1}^*$) be the standard basis vectors for the first (resp.\@ second) copy of $\C^{n-1}$ in the $\SO(n-1,\C)$-decomposition\footnote{Self-duality of the $\SO(n-1,\C)$-representation $\C^{n-1}$ with respect to the standard Euclidean bilinear form ensures that this notation is compatible with Theorem \ref{fft_theorem}.}, the $i$\textsuperscript{th} copy of $\C^2$ in the $\Sp(2,\C)$-decomposition (which we shall henceforth denote by $(\C^2)_i$) has basis $e_i,e_i^*$.
	%
	%
	We now consider separately the two subcases $\Lambda_{\inv}^{k}U$ with $k=n-1$ and $k\neq n-1$. It suffices to show that (\ref{ExteriorUdimensions}) holds in both subcases, as the result will then follow by the same argument as in Case 2.
	\begin{enumerate}[(i)]
		\item \textul{$k=n-1$:} To begin, we note that the symplectic form on $(\C^2)_i$ stabilized by $\Sp(2,\C)$ is given by
		\begin{align}\label{pisymplecticformdefinition}
			p_i:= e_i\otimes e_i^* - e_i^* \otimes e_i.
		\end{align}
		For each pair $i,j\in\{1,\dots,n-1\}$, the copies $(\C^2)_i$ and $(\C^2)_j$ are isomorphic as $\Sp(2,\C)$-modules via $f_{i,j}\: e_i\mapsto e_j, e_i^* \mapsto e_j^*$, and the natural extension of this map to tensors satisfies $f_{i,j}(p_i)=p_j$. By Theorem \ref{fft_theorem}(ii), the $\Sp(2,\C)$-invariant tensors in $T(\C^2)$ are spanned by mutations of tensor powers of the symplectic form defining the group. The tensor algebra of $U=\bigoplus_{i=1}^{n-1} (\C^2)_i$ is the direct sum of all possible tensor products of the spaces $(\C^2)_i$, $i=1,\dots n-1$, hence the $\Sp(2,\C)$-invariant tensors are spanned by mutations of tensor products of the symplectic forms $p_i$, $i=1,\dots, n-1$. The symplectic forms $p_i$ have degree $2$, hence $(U^{\otimes (n-1)})_{\inv}$ can only be non-trivial if $n-1$ is even. We assume for the rest of the subcase that $n-1=2l$. Multilinearly expanding a tensor product of the form $p_{i_1}\otimes \dots \otimes p_{i_{l}}$, or any mutation thereof, one sees that the vectors $e_{i_s}, e_{i_s}^*$, $s=1,\dots l$ appear in each term. Similarly, each term of the $\SO(n-1,\C)$-invariant tensor $I=\sum_{i=1}^{n-1} e_i\otimes e_i^*$ from Theorem \ref{fft_theorem}(i) contains a pair of vectors of the form $e_i,e_i^*$. The flips of $I$ are precisely
		\begin{align}\label{Chap_Hof22:flips_of_I}
			I=\sum_{i=1}^{n-1} e_i\otimes e_i^*, \quad I_1 := \sum_{i=1}^{n-1} e_i^*\otimes e_i, \quad I_2:= \sum_{i=1}^{n-1} e_i\otimes e_i, \quad I_3:=  \sum_{i=1}^{n-1} e_i^*\otimes e_i^*,
		\end{align}hence any mutation of an $l$-fold tensor product of flips of $I$ has the property that, when fully expanded, each term contains a pair of vectors of the form $e_i,e_i^*$ or $e_i,e_i$ or $e_i^*, e_i^*$ for some $i=1,\dots, n-1$. On the other hand, the $\SO(n-1,\C)$-invariant tensor \[ \det= \sum_{\sigma\in S_{n-1}} \text{sign}(\sigma) e_{\sigma(1)}\otimes \dots \otimes e_{\sigma(n-1)}\] from Theorem \ref{fft_theorem}(i) (and any mutation and/or flip thereof) has the property that, for each $i=1,\dots, n-1$, each term contains either exactly one copy of $e_i$ or exactly one copy of $e_i^*$. Said differently, linear combinations of mutations of $l$-fold tensor products of $p_i,I,I_1,I_2,I_3$ have repeated indices in each term, whereas linear combinations of mutations of flips of $\det$ do not have any terms with repeated indices. Comparing Theorem \ref{fft_theorem}(i) with the above observation that the $\Sp(2,\C)$-invariant tensors are spanned by mutations of tensor products of the $p_i$, we then have that $(U^{\otimes (n-1)})_{\inv}$ is contained in the span of all mutations of ($l$-fold) tensor products of flips of $I$. Consider the skew-symmetrization map $\Alt\: T(U) \to \Lambda^{\bullet}U \subset T(U)$ given by 
		\[
		u_{i_1}\otimes \dots \otimes u_{i_k}\mapsto \sum_{\sigma \in S_k} \text{sign}(\sigma) u_{\sigma(i_1)} \otimes \dots \otimes u_{\sigma(i_k)}.
		\]
		This map is clearly $H$-equivariant and maps any invariant exterior form (viewed as a skew-symmetric tensor) to a positive multiple of itself, so it suffices to consider the image of $(U^{\otimes 2l})_{\inv}$ under $\Alt$. Noting that the skew-symmetrization of a tensor and the skew-symmetrization of any mutation of that tensor agree up to sign, it follows that $\Lambda_{\inv}^{2l}U$ is contained in the span of all images under $\Alt$ of $l$-fold tensor products of flips of $I$. Furthermore, since for any tensors $\alpha,\beta$ the skew-symmetric tensor $\Alt(\alpha\otimes \beta)$ agrees up to positive scaling with $\Alt(\alpha)\wedge \Alt(\beta)$, we have
		\[
		\Lambda_{\inv}^{2l}U\subseteq \Lambda^l \mathcal{S},\qquad \text { where }   \mathcal{S}:=\Span_{\C}\{\Alt(I),\Alt(I_1),\Alt(I_2),\Alt(I_3)\}.
		\]One easily sees that $\Alt(I) = I-I_1= - \Alt(I_1)$, $\Alt(I_2) = \Alt(I_3)=0$, and we note furthermore that the tensor 
		\[
		\mathcal{I}:=\Alt(I) =  \sum_{i=1}^{n-1} (e_i\otimes e_i^* - e_i^*\otimes e_i) = \sum_{i=1}^{n-1} p_i
		\]is $\Sp(2,\C)$-invariant. Thus $\mathcal{S}= \C \mathcal{I}$ is the trivial $H^{\C}$-representation, and it follows that
		 \[\Lambda_{\inv}^{2l}U= \Lambda^l \mathcal{S} = \Span_{\C}\{ \underbrace{\mathcal{I}\wedge \dots \wedge \mathcal{I}}_{\text{($l$ copies)}} \}.\]
		In particular we have shown that $\dim_{\C}\Lambda_{\inv}^{n-1}U $ is equal to $1$ if $n-1$ is even and $0$ if $n-1$ is odd, as desired.
		\item \textul{$k\neq n-1$:} Similarly to (\ref{Case2eqn1}), taking the $k$\textsuperscript{th} exterior power of the $\SO(n-1,\C)$-decomposition $U\simeq \C^{n-1}\otimes \C^{n-1}$ gives
		\begin{align}\label{Case3eqn1}
			\Lambda^k U \simeq  \Lambda^k( \C^{n-1} \oplus \C^{n-1}) \simeq \bigoplus_{p+q=k}(\Lambda^p \C^{n-1})\otimes (\Lambda^q \C^{n-1})
		\end{align}as $\SO(n-1,\C)$-modules. Recalling from \cite[Eqn.\@ (24)]{homdata} that $\left( \begin{smallmatrix}
		1&0 \\ 0 &-1
	\end{smallmatrix}\right)  \in \mathfrak{h}^{\C} $ acts on $U$ by $e_i\mapsto e_i$, $e_i^*\mapsto -e_i^*$, it follows that the invariant elements in (\ref{Case3eqn1}) must lie in a summand with $p=q$ (which is only possible when $k$ is even). The possibility $p=q=\frac{n-1}{2}$ is excluded by the assumption $k\neq n-1$, and we note that if $p=q\neq \frac{n-1}{2}$ then the $\SO(n-1,\C)$-module $\Lambda^p \C^{n-1}\simeq  \Lambda^q\C^{n-1}$ is irreducible (see Theorem \ref{GWinvarianttheorem}(i)) and self-dual. In particular, for $k\neq n-1$ the dimension of the space of invariant $k$-forms on $U$ satisfies the upper bound
	\begin{align}\label{ExteriorUdimensions_inequality}
	\dim_{\C} \Lambda^k_{\inv} (U) \leq  \begin{cases*} 1 & \text{if $k$ is even and $k\leq \dim_{\C} U$,} \\ 0 & \text{otherwise.} \end{cases*}
\end{align}
	But the form $\omega_{1,0}$ (resp.\@ $\omega_{0,1}$) considered in Case $2$ is again an element of $\Lambda_{\inv}^{2,0}(\mathfrak{m}^{\C}_{\mathcal{H}})$ (resp.\@ $\Lambda_{\inv}^{0,2}(\mathfrak{m}^{\C}_{\mathcal{H}})$) whose top power is a $\varphi_1$-holomorphic (resp.\@ $\varphi_1$-anti-holomorphic) volume form, hence equality is achieved in (\ref{ExteriorUdimensions_inequality}), as desired.
	\end{enumerate}

	\textul{Case 4: The Five Exceptional Spaces.}
	These are the five spaces from Theorem \ref{hom3Sas_classification} with $G$ an exceptional Lie group. Following \cite{homdata}, we denote the corresponding $\ts$ data by
	\[
	\mathfrak{g}^s=\mathfrak{g}_0^s,\oplus \mathfrak{g}_1^s,\qquad (\mathfrak{g}_1^s)^{\C} \cong \C^2\otimes U^s, \qquad s=1,2,3,4,5,
	\]
	and we recall that the $(\mathfrak{h}^s)^{\C}$-modules $U^s$ are given as highest weight modules on \cite[p.\@ 841]{homdata}. This information is summarized in Table \ref{Chap_Hof22:Tab:exceptionalcasesUrep}.
	\begin{table}[h!] 
		\centering
		\begin{tabular}{ |l||l|l|l|l|l| }
			\hline 
			&  $s=1$ & $s=2$ & $s=3$ & $s=4$ & $s=5$ \\
			\hline \hline 
			$G^s$   &    $\G_2$ & $\text{F}_4$ & $\text{E}_6$ & $\text{E}_7$ & $\text{E}_8$ \\ 
			\hline
			$H^s$   &    $\Sp(1)$ & $\Sp(3)$ & $\SU(6)$ & $\Spin(12)$ & $\text{E}_7$ \\ \hline
			$(\mathfrak{h}^s)^{\C}$   &    $A_1=\mathfrak{sp}(2,\C)$ & $C_3=\mathfrak{sp}(6,\C)$ & $A_5 = \mathfrak{sl}(6,\C)$ & $D_6=\mathfrak{so}(12,\C)$ & $\text{E}_7= \mathfrak{e}_7^{\C}$ \\ \hline
			$U^s$   &    $V(3)$ & $V(\lambda_3)$ & $V(\lambda_3)$ & $V(\lambda_5)$ & $V(\lambda_7)$\\
			\hline
		\end{tabular}
		\caption{The Exceptional Homogeneous $\ts$ Spaces}
		\label{Chap_Hof22:Tab:exceptionalcasesUrep}
	\end{table}
	Using the LiE computer algebra package (\cite{LiE}), one checks that (\ref{ExteriorUdimensions}) holds for each $s=1,2,3,4,5$,
	and the result in this case then follows by the same argument as in Case 2.
\end{proof}
As a consequence, we immediately obtain a description of the invariant $\varphi_1$-(anti-)holomorphic\footnote{Clearly $\varphi_1$ is not an almost complex structure since the dimension of the manifold is odd; by \emph{$\varphi_1$-(anti-)holomorphic tangent spaces} we mean the $\pm i$ eigenspaces of $\varphi_1$ on the complexification of $(\R\xi_1)^{\perp}$, and \emph{$\varphi_1$-(anti-)holomorphic forms} may then be defined relative to this subspace in the usual way.} forms on the full tangent bundle:
\begin{theorem}\label{holomorphictheorem}
	If $(M=G/H, g, \xi_i,\eta_i,\varphi_i)$ is a simply-connected homogeneous $\ts$ space, then the invariant $\varphi_1$-(anti-)holomorphic forms are given by
	\begin{align*}
		\Lambda_{\inv}^{\bullet,0 }(\mathfrak{m}^{\C})= \Span_{\C}\{\omega_{1,0}^k, \ y_{1,0}\wedge \omega_{1,0}^k \}_{k=0}^{n-1} , \qquad  
		\Lambda_{\inv}^{0, \bullet }(\mathfrak{m}^{\C})=\Span_{\C}\{\omega_{0,1}^k, \ y_{0,1}\wedge \omega_{0,1}^k \}_{k=0}^{n-1},
	\end{align*}
	where 
	\[
	y_{1,0}:= (\xi_2-i\xi_3), \quad  y_{0,1}:=\overline{y_{1,0}}, \qquad \qquad   \omega_{1,0}:= (\Phi_2\rvert_{\mathcal{H}}-i\Phi_3\rvert_{\mathcal{H}}), \quad  \omega_{0,1}:= \overline{\omega_{1,0}}.
	\]
\end{theorem}
\begin{proof}
	Since the isotropy group acts trivially in the vertical directions, we have
	\begin{align*}
		\Lambda_{\inv}^{\bullet,0}(\mathfrak{m}^{\C}) &= \Lambda^{\bullet,0}(\mathfrak{m}_{\mathcal{V}}^{\C})\otimes \Lambda_{\inv}^{\bullet,0}(\mathfrak{m}_{\mathcal{H}}^{\C}) \cong [ 1\otimes  \Lambda_{\inv}^{\bullet,0}(\mathfrak{m}_{\mathcal{H}}^{\C}) ] \oplus[ y_{1,0}\otimes \Lambda_{\inv}^{\bullet,0}(\mathfrak{m}_{\mathcal{H}}^{\C})] , \\
		\Lambda_{\inv}^{0,\bullet}(\mathfrak{m}^{\C}) &= \Lambda^{0,\bullet}(\mathfrak{m}_{\mathcal{V}}^{\C})\otimes \Lambda_{\inv}^{0,\bullet}(\mathfrak{m}_{\mathcal{H}}^{\C}) \cong [ 1\otimes  \Lambda_{\inv}^{0,\bullet}(\mathfrak{m}_{\mathcal{H}}^{\C}) ] \oplus[ y_{0,1}\otimes \Lambda_{\inv}^{0,\bullet}(\mathfrak{m}_{\mathcal{H}}^{\C})] ,
	\end{align*}
	and the result then follows from an analogous argument as in Case $2$ in the proof of the preceding proposition, where it was noted that $ \omega_{1,0}$ (resp.\@ $ \omega_{0,1}$) generates $\Lambda_{\inv}^{\bullet,0}(\mathfrak{m}_{\mathcal{H}}^{\C})$ (resp.\@ $\Lambda_{\inv}^{0,\bullet}(\mathfrak{m}_{\mathcal{H}}^{\C})$).
\end{proof}

Using Theorem \ref{holomorphictheorem}, we are now ready to prove the main result of the section:
\begin{theorem}\label{invspinortheorem}
	For a simply-connected homogeneous $3$-Sasakian manifold $(M=G/H,g,\xi_i,\eta_i,\varphi_i)$ of dimension $4n-1$, the invariant spinors are given with respect to any adapted basis by 
	\begin{align}\label{invspinorbasis}
	 \Sigma_{\inv} = \Span_{\C}\{ \omega^k, \ y_1\wedge \omega^k\}_{k=0}^{n-1},
	 \end{align} where $\omega := \sum_{p=1}^{n-1}y_{2p}\wedge y_{2p+1}$. 
\end{theorem}
\begin{proof}
	It is well-known that the $\ts$ structure on $M$ gives a reduction of the structure group of the tangent bundle to $\Sp(n-1) \subset \SO(4n-1)$ (see \cite[Thm.\@ 5]{3Sas_structure_reduction}). Furthermore, since the $\ts$ structure on $M=G/H$ is assumed to be homogeneous, we have that the image of the isotropy representation $\iota$ is contained in this reduction: 
	\begin{align}\label{spn_reduction_isotropy_containment}
	\iota(H)\subseteq \Sp(n-1) \subseteq \SO(4n-1).
	\end{align}The lifted action of $\Sp(n-1)$ on $\Sigma=\Lambda^{\bullet}L'$ is given explicitly in \cite[Prop.\@ 4.5]{AHLspheres} (see also the discussion in \cite[Section 4.1.4]{AHLspheres}): since operators in $\mathfrak{sp}(2n-2,\C)$ are traceless, one has $\Sigma\simeq \Lambda^{\bullet}L'$ as $\Sp(n-1)$-representations. The containment (\ref{spn_reduction_isotropy_containment}) then implies that $\Sigma\simeq \Lambda^{\bullet}L'$ also holds as $H$-representations, hence $ \Sigma M \cong \Lambda^{0,\bullet}(T^*_{\C}M) $ as homogeneous vector bundles, and the result follows from Theorem \ref{holomorphictheorem} by noting that $y_1=\frac{1}{\sqrt{2}}y_{0,1}$ and $\omega = -\frac{1}{2}\omega_{0,1}$.  
\end{proof}

\begin{remark}
	The motivation for Theorem \ref{invspinortheorem} can be found in \cite[Eqn.\@ (37), Thm.\@ 4.11]{AHLspheres}, which prove the result for the case $S^{4n-1}=\Sp(n)/\Sp(n-1)$. Indeed, these spheres are in many ways the model examples of homogeneous $\ts$ geometries, so it is not surprising that the result holds in general.
\end{remark}

In the proof of Theorem \ref{invspinortheorem} we relied on the identification $\Sigma \cong \Lambda^{0,\bullet}_{\C}(M) $  between invariant $\varphi_1$-anti-holomorphic forms and invariant spinors. A natural question then arises as to the relationship between invariant real differential forms and invariant spinors. The next few results are devoted to the exploration of this topic. From this point forward we fix, without further mention, the Clifford algebra representation associated to an adapted basis, so that the invariant spinors take the form (\ref{invspinorbasis}).
\begin{remark} Recall from \cite[Section 4.5]{homdata} that the space $\SU(n+1)/ \text{S}(\U(n-1)\times \U(1))$ has an additional invariant $2$-form compared to the other homogeneous $\ts$ spaces, which is given in an adapted basis by \[\Phi_0 = \sum_{p=1}^{n-1} (e_{4p}\wedge e_{4p+1} - e_{4p+2}\wedge e_{4p+3}).\]
\end{remark}
\begin{lemma}\label{Chap_Hof22:cliffproductinvforms}Let $(M=G/H,g,\xi_i,\eta_i,\varphi_i)$ be a simply-connected homogeneous $3$-Sasakian manifold of dimension $4n-1$. Then for any integer $k\geq 0$ we have
	\begin{align}
		\Phi_0\cdot \omega^k &= 0 , \\
		(\Phi_1\rvert_{\mathcal{H}})\cdot \omega^k &= 2i(2k-n+1) \omega^k , \\
		\label{cliffproductPhi2} (\Phi_2\rvert_{\mathcal{H}})\cdot \omega^k &=  2(\omega^{k+1} -k (n-k) \omega^{k-1}), \\
		(\Phi_3\rvert_{\mathcal{H}})\cdot \omega^k &=  -2i(\omega^{k+1} +k (n-k) \omega^{k-1}),
	\end{align}
	where $\omega := \sum_{p=1}^{n-1} y_{2p}\wedge y_{2p+1}$ and we use the convention $\omega^0=1$.
\end{lemma}
\begin{proof}
			%
			%
			%
			%
			%
			%
			
			These identities follow from a straightforward calculation in the spin representation. Considering first $\Phi_0$, for $k\geq 0$ we calculate:
			\begin{align*}
				\Phi_0 \cdot \omega^k &= \sum_{p=1}^{n-1} (e_{4p}\wedge e_{4p+1} - e_{4p+2}\wedge e_{4p+3}) \cdot \omega^k  \\
				&= i\sum_{p=1}^{n-1} [(x_{2p}\lrcorner + y_{2p}\wedge )( y_{2p}\wedge -x_{2p}\lrcorner)  - (  x_{2p+1}\lrcorner+y_{2p+1}\wedge)(y_{2p+1}\wedge -  x_{2p+1}\lrcorner) ]  \omega^k \\
				&= i\sum_{p=1}^{n-1} [ x_{2p}\lrcorner(y_{2p}\wedge \omega^k ) -y_{2p}\wedge (x_{2p}\lrcorner\omega^k)   -x_{2p+1}\lrcorner(y_{2p+1}\wedge \omega^k) + y_{2p+1}\wedge (x_{2p+1}\lrcorner\omega^k)     ]     \\
				&=  i\sum_{p=1}^{n-1} [  \omega^k - 2ky_{2p}\wedge y_{2p+1}\wedge \omega^{k-1}    - \omega^k + 2k y_{2p+1} \wedge (-y_{2p}) ]  \\
				&= 0.
			\end{align*}
			The calculations for $(\Phi_i\rvert_{\mathcal{H}})\cdot \omega^k$, $i=1,2,3$ are analogous.
		\end{proof}
We immediately deduce:
\begin{corollary}\label{S0S1S2S3corr}
	For $i\in\{0,1,2,3\} $, let $S_i$ denote the complex span of the spinors $(\Phi_i\rvert_{\mathcal{H}})^k \cdot 1 $ with $k= 1,\dots, 2n-1$. We have:
	\begin{align*}
		S_0=\{0\}, \quad S_1 = \Span_{\C}\{1\}, \quad S_2=S_3 = \Span_{\C}\{\omega^k\}_{k=0}^{n-1}.
	\end{align*} 
\end{corollary}
\begin{proof}
	The cases $S_0$ and $S_1$ are clear from the preceding lemma. From (\ref{cliffproductPhi2}) we note that the Clifford product of the form $\frac{1}{2}\Phi_2\rvert_{\mathcal{H}} $ with the spinor $\omega^k$ is a monic degree $(k+1)$ polynomial in $\omega$, hence $\frac{1}{2^k}(\Phi_2\rvert_{\mathcal{H}})^k \cdot 1$ is a monic degree $k$ polynomial in $\omega$. By induction, one sees that
	\[
	\Span_{\C}\{ \tfrac{1}{2^k}(\Phi_2\rvert_{\mathcal{H}})^k\cdot 1 \}_{k=0}^{k_0} = \Span_{\C}\{ \omega^k\}_{k=0}^{k_0}
	\]for any $k_0\in \{1,\dots, n-1\}$, and the result for $S_2$ follows. The argument for $S_3$ is analogous.
\end{proof}

This also gives a nice description of the invariant spinors in terms of the invariant real differential forms:
\begin{theorem}
	The space $\Sigma_{\inv}$ of invariant spinors on a simply-connected homogeneous $\ts$ manifold is spanned by Clifford products of invariant differential forms with the invariant spinor $1\in \Sigma_{\inv}$. 
\end{theorem}
\begin{proof}
	In light of Theorem \ref{invspinortheorem}, it suffices to show that spinors of the form $\omega^k$ and $y_1\wedge \omega^k$ can be obtained as linear combinations of Clifford products of invariant differential forms with $1\in \Sigma_{\inv}$. This follows from Corollary \ref{S0S1S2S3corr} by noting that $\omega^k\in S_2$ and $y_1\wedge \omega^k \in \xi_2\cdot S_2$. 
\end{proof}
\begin{remark}
	We would like to point out that the results of this section so far also hold in the more general setting of compact simply-connected homogeneous $\tad$ spaces; The reason for this is that the generalized $3$-Sasakian data used to define homogeneous $\tad$ structures coincides, in the case of a compact space, with the notion of $\ts$ data (compare Theorem \ref{DOP_3Sas_data} with \cite[Section 3.1]{hom3alphadelta}. In particular, the isotropy representation of a family of compact homogeneous $\tad$ spaces parameterized by $\alpha,\delta>0$ is isomorphic to the isotropy representation of the corresponding homogeneous $\ts$ space obtained by setting $\alpha=\delta=1$. The next section discusses Killing spinors on homogeneous $\ts$ spaces, which do not carry over to the corresponding $\tad$ spaces. Rather, certain deformations of Killing spinors in the $\tad$ setting are investigated in \cite{AHduality}. 
\end{remark}

\section{The Space of Riemannian Killing Spinors}\label{KSsection}

We conclude the paper with an explicit description of the Riemannian Killing spinors on a simply-connected homogeneous $\ts$ space:
\begin{theorem}\label{KSbasis_theorem}
	Let $(M^{4n-1}=G/H,$ $ g,\xi_i,\eta_i,\varphi_i)$ be a simply-connected homogeneous $\ts$ manifold, and fix a description of the spinor module relative to an adapted basis as in the previous section. If $n\geq 2$, then the space of invariant Killing spinors has a basis given by
	\begin{align}
\label{invKSbasis}	\psi_{k}:= \omega^{k+1}-i(k+1)y_1\wedge \omega^k, \qquad  -1\leq k \leq n-1,
	\end{align}
	where we use the conventions $\omega^{-1}=0$ and $\omega^0=1$. If $n=1$ then the space of invariant Killing spinors has a basis given by $1$, $y_1$. Furthermore, if $(M,g)\ncong (S^{4n-1},g_{\round})$ then any Killing spinor is invariant.
\end{theorem}
\begin{proof}
	Let $\Uplambda$, $\Uplambda^g \: \mathfrak{m}\to \mathfrak{so}(\mathfrak{m})$ denote the Nomizu maps of the canonical and Levi-Civita connections respectively, and note that the Nomizu maps $\widetilde{\Uplambda}$, $\widetilde{\Uplambda^g}$ of the corresponding spinorial connections are given by composing with the isomorphism $\mathfrak{spin}(\mathfrak{m})\cong \mathfrak{so}(\mathfrak{m})$. First we consider the horizontal directions $X\in\mathcal{H}$. From (\ref{canonicalNomizumap}) we have $\Uplambda(X)=0$ for all $X\in\mathcal{H}$, and thus any invariant Killing spinor $\psi$ satisfies the algebraic equation
	\begin{align}\label{RKSeqn}
		0= \widetilde{\Uplambda}(X)\cdot \psi = \widetilde{\Uplambda^g}(X)\cdot \psi +\frac{1}{4}(X\lrcorner T) \cdot \psi = \frac{1}{2}X\cdot \psi + \frac{1}{2} \sum_{i=1}^3 \xi_i\cdot \varphi_i(X)\cdot \psi \qquad 	\text{for all } X\in \mathcal{H}.
	\end{align}
Assume first that $n\geq 2$. Calculating in an adapted basis $\{e_{4p},e_{4p+1},e_{4p+2},e_{4p+3}\}$, we find:
		\begin{align*}
			\frac{1}{2} e_{4p}\cdot \omega^k +& \frac{1}{2}\sum_{i=1}^3 \xi_i\cdot \varphi_i(e_{4p})\cdot \omega^k = 	\frac{1}{2}  \left( e_{4p}\cdot \omega^k + \xi_1\cdot e_{4p+1}\cdot \omega^k + \xi_2\cdot e_{4p+2}\cdot \omega^k  +\xi_3\cdot e_{4p+3}\cdot \omega^k     \right)\\
			&=\frac{1}{2}  \left( 2i  x_{2p} \lrcorner\omega^k + \xi_2\cdot e_{4p+2}\cdot \omega^k  +\xi_3\cdot e_{4p+3}\cdot \omega^k     \right)\\
			&= i  x_{2p}\lrcorner \omega^k +\frac{1}{2}\left( -y_1\wedge (y_{2p+1}\wedge +x_{2p+1}\lrcorner  )\omega^k +  y_1\wedge (y_{2p+1}\wedge - x_{2p+1}\lrcorner  )\omega^k      \right) \\
			&= i x_{2p}\lrcorner\omega^k - y_1\wedge (x_{2p+1} \lrcorner \omega^k) ,
		\end{align*}
	and similarly, {\small
		\begin{align*}
			\frac{1}{2}e_{4p}\cdot (y_1\wedge \omega^k) + \frac{1}{2} \sum_{i=1}^3 \xi_i\cdot \varphi_i(e_{4p})\cdot (y_1\wedge \omega^k) 
			&= iy_{2p}\wedge (y_1\wedge \omega^k) +y_{2p+1}\wedge \omega^k .
		\end{align*}
	}Writing
	\[\psi:= \sum_{k=0}^{n-1}\lambda_k \omega^k + \sum_{k=0}^{n-1} \lambda'_k (y_1\wedge \omega^k)
	\]
	in terms of the basis from Theorem \ref{invspinortheorem}, and assuming that $\psi$ is a Killing spinor, we have
		\begin{align}
			0&= \widetilde{\Uplambda}(e_{4p})\cdot \psi =  \frac{1}{2}e_{4p}\cdot \psi +\frac{1}{2}\sum_{i=1}^3 \xi_i\cdot \varphi_i(e_{4p})\cdot \psi \nonumber \\
		\label{KSmastereqn}	&= \sum_{k=0}^{n-1} \lambda_k[ix_{2p}\lrcorner\omega^k -y_1\wedge(x_{2p+1}\lrcorner\omega^k)  ]  + \sum_{k=0}^{n-1} \lambda_k'[y_{2p+1}\wedge \omega^k+  iy_{2p}\wedge(y_1\wedge \omega^k) ] . 
		\end{align}
	%
	For the $k=l$ index of the summations in (\ref{KSmastereqn}), the degrees of the four terms are $2l-1$, $2l$, $2l+1$ and $2l+2$ respectively. Considering separately the even and odd degree parts of (\ref{KSmastereqn}), we are seeking solutions of
	\begin{align*}
		\sum_{k=0}^{n-1} [ i\lambda_k  x_{2p}\lrcorner\omega^k + \lambda_k' y_{2p+1}\wedge \omega^k ] = 0 = \sum_{k=0}^{n-1} [ -\lambda_k y_1\wedge (x_{2p+1}\lrcorner \omega^k) +i\lambda_k' y_{2p}\wedge (y_1\wedge \omega^k)   ],
	\end{align*}
	or equivalently, solutions of the linear system of equations
	\[
	\lambda_{k}' = -i(k+1)\lambda_{k+1}, \qquad -1\leq k \leq n-1.
	\]
	This gives $n+1$ linearly independent spinors
	\[
	\psi_{k}:=  \omega^{k+1} -i(k+1) y_1\wedge \omega^k , \qquad -1\leq k \leq n-1,
	\]
	and a straightforward calculation of the other horizontal derivatives (by substituting $X=e_{4p+1}, e_{4p+2}$, and $e_{4p+3}$ into (\ref{RKSeqn})) shows that these spinors satisfy the Killing equation in the horizontal directions. In the vertical directions, one sees from Theorem \ref{DOP_3Sas_data} that the Nomizu map for the Levi-Civita connection satisfies $\Uplambda^g(\xi_i) = \xi_j\wedge \xi_k$ for any even permutation $(i,j,k)$ of $(1,2,3)$. Considering the spin lifts $\widetilde{\Uplambda^g}(\xi_i)=\frac{1}{2}\xi_j\cdot \xi_k$, another straightforward calculation in the spin representation shows that any spinor of the form $\omega^k$ or $y_1\wedge \omega^k$ satisfies the Killing equation in the vertical directions, and we conclude that the $\psi_{k}$ are Killing spinors in the case $n\geq 2$. Assuming now that $n=1$, the horizontal distribution is trivial, and consequently Equation (\ref{RKSeqn}) does not apply. In this dimension the spinor bundle has complex dimension equal to $2$, hence it is spanned by $1, y_1$ and the above argument for the vertical directions shows that these are Killing spinors. The dimension of the space of Killing spinors on a $\ts$ manifold $(M^{4n-1}, g)\ncong (S^{4n-1},g_{\text{round}})$ is equal to $n+1$ (see \cite{Bar}), and the result follows.
\end{proof}
\begin{remark}
	The final assertion of Theorem \ref{KSbasis_theorem}--that any Killing spinor on a homogeneous $\ts$ space $(M^{4n-1}=G/H, g)\ncong (S^{4n-1},g_{\text{round}})$ is invariant--was previously proved using a different method in \cite[Thm.\@ 7.1]{kath_Tduals}: Kath showed that $G$ has a representation on the space of Killing spinors, and then deduced that any Killing spinor is invariant (equivalently, this representation is trivial) by comparing in each case the dimension of the space of Killing spinors with the dimension of the smallest non-trivial representation of $G$. 
\end{remark}
\begin{remark}
	Let us comment briefly on the case of the round sphere. Using B\"{a}r's correspondence in \cite{Bar} between Killing spinors on a simply-connected manifold and parallel spinors on its cone, it is easy to see that the spinor bundle of the round sphere is parallelized by Killing spinors; the cone over the round sphere is the punctured Euclidean space of one dimension higher, which has trivial holonomy and therefore a parallelization of its spinor bundle by parallel spinors. In particular, for the $\ts$ round sphere $(S^{4n-1}=\frac{\Sp(n)}{\Sp(n-1)},g_{\round})$, the $\Sp(n)$-invariant Killing spinors $\psi_k$ from Theorem \ref{KSbasis_theorem} fail to span the whole spinor bundle when $n>1$, so there exist non-invariant Killing spinors in this case. An explicit formula for the Killing spinors on the round sphere in stereographic coordinates can be found in \cite[Ex.\@ 2 on p.\@ 37]{BFGK}.
\end{remark}
We conclude the section by exploring which of the invariant Killing spinors from Theorem \ref{KSbasis_theorem} recover the invariant $\ts$ structure on $(M=G/H,g,\xi_i,\eta_i,\varphi_i)$ via the construction in Section \ref{constructionsection}. In light of Theorem \ref{3Sasakian_converse} and Definition \ref{associated_tensors_definition}, we see that in order to recover the Sasakian structure $(\xi_i,\eta_i,\varphi_i)$ it suffices to find a Killing spinor $\Psi$ such that $\Psi':= -\xi_i\cdot \Psi$ is also a Killing spinor. Thus, in the homogeneous $\ts$ setting, we consider the subbundles $\mathcal{E}_i$ of the spinor bundle spanned (over the space $C^{\infty}(M)$ of smooth real-valued functions) by invariant spinors with this property:
\[
\mathcal{E}_i:= \Span_{C^{\infty}(M)} \{ \Psi \in \kappa_{\inv}(M,g) \: \ \Psi':= -\xi_i\cdot \Psi \in \kappa_{\inv}(M,g)\},\qquad  i=1,2,3,
\]
where $\kappa_{\inv}(M,g):=\Span_{\C}\{\psi_k\}_{k=-1}^{n-1}$ denotes the space of invariant Killing spinors. In fact, it turns out that these subbundles coincide with the bundles $E_i^-$ from (\ref{Eplusminusspaces_3Sas}), as we prove in the following proposition:
\begin{proposition}\label{Chap_Hof22:E_spaces_bases_hom_ts}
	If $(M^{4n-1}=G/H,g,\xi_i,\eta_i,\varphi_i)$ is a simply-connected homogeneous $\ts$ space, then $\mathcal{E}_i= E_i^-$ for $i=1,2,3$. Furthermore, each $\mathcal{E}_i$ has a basis $\Psi_{\mathcal{E}_i,0},\Psi_{\mathcal{E}_i,1}$ given by
	\begin{align*}
		\Psi_{\mathcal{E}_1,0}&:= \psi_{-1}, \quad \Psi_{\mathcal{E}_1,1}:= \frac{i}{n} \psi_{n-1}, \quad \Psi_{\mathcal{E}_2,0} := \sum_{k=0}^{\lfloor \frac{n-1}{2} \rfloor} \frac{(-1)^k}{(2k+1)!} \psi_{2k},\\
		\Psi_{\mathcal{E}_2,1}&:= \sum_{k=0}^{\lfloor \frac{n}{2} \rfloor} \frac{(-1)^k}{(2k)!} \psi_{2k-1},  \quad \Psi_{\mathcal{E}_3,0} := \sum_{k=0}^{\lfloor \frac{n-1}{2} \rfloor} \frac{1}{(2k+1)!} \psi_{2k}, \quad  \Psi_{\mathcal{E}_3,1}:= \sum_{k=0}^{\lfloor \frac{n}{2} \rfloor} \frac{1}{(2k)!} \psi_{2k-1}, 
	\end{align*}
where for $n=1$ we adopt the conventions $\psi_{-1} = 1$, $\psi_0 = y_1$.
\end{proposition}
\begin{proof}
	Suppose that $n\geq 2$. Considering first the case $(M^{4n-1}, g)\ncong (S^{4n-1},g_{\text{round}})$, Theorem \ref{KSbasis_theorem} implies that any Killing spinor is invariant, and it then follows from the proof of Theorem \ref{Sasakian_converse} that $E_i^-\subseteq \mathcal{E}_i$. Therefore, in order to show that $E_i^-=\mathcal{E}_i$ it suffices to show that $\text{rank}(\mathcal{E}_i) = 2$. 
	To find basis elements $\Psi\in\mathcal{E}_i$, we write $\Psi = \sum_{k=-1}^{n-1} \lambda_k \psi_k$ in terms of the basis (\ref{invKSbasis}) of invariant Killing spinors and we seek to determine for which values of $\lambda_{-1},\dots, \lambda_{n-1}$ there exist $\Theta_{-1},\dots,\Theta_{n-1}\in\C$ satisfying
	\begin{align}\label{inducingKSMasterEqn}
		\xi_i \cdot \Psi = \sum_{k=-1}^{n-1} \Theta_k \psi_k. 
	\end{align}
	To show that $\text{rank}(\mathcal{E}_i)= 2$ for $i=1,2,3$, we treat the subcases $i=1,2,3$ individually:
	\begin{enumerate}[(i)]
		\item \textul{$i=1$:} First, we note that 
		\[
		\xi_1\cdot \psi_k = i\omega^{k+1} - (k+1) y_1\wedge \omega^k.
		\]
		Substituting this into (\ref{inducingKSMasterEqn}), we are looking for solutions of
		\[
		\sum_{k=-1}^{n-1} \lambda_k[i\omega^{k+1} - (k+1)y_1\wedge \omega^k] =  \sum_{k=-1}^{n-1} \Theta_k[\omega^{k+1}-i(k+1)y_1\wedge \omega^k] ,
		\] or equivalently, solutions of the linear equations 
		\[
		i\lambda_{-1} = \Theta_{-1} , \quad \lambda_{n-1} = i \Theta_{n-1}, \quad i\lambda_k = \Theta_k, \quad \lambda_k = i\Theta_k, \quad \text{ for }k=0,\dots, n-2.
		\]
		The solutions of this system of equations necessarily have $\lambda_k=\Theta_k=0$ for all $k\neq -1,n-1$ (and no restriction on the other variables), hence are spanned by $\Psi_{\mathcal{E}_1,0}= 1$ and $\Psi_{\mathcal{E}_1,1}=y_1\wedge \omega^{n-1}$.
		\item \textul{$i=2$:} Proceeding similarly as in the previous subcase, we first note that
		\[
		\xi_2\cdot \psi_k = iy_1\wedge \omega^{k+1} +(k+1)\omega^k.
		\]
		Substituting this into (\ref{inducingKSMasterEqn}) gives
		\[
		\sum_{k=-1}^{n-1} \lambda_k[iy_1\wedge \omega^{k+1} +(k+1)\omega^k] =   \sum_{k=-1}^{n-1} \Theta_k[\omega^{k+1}-i(k+1)y_1\wedge \omega^k],
		\]
		or equivalently, the linear system
		\[
		(k+1)\lambda_k=\Theta_{k-1},\quad \lambda_{k-1}=-(k+1)\Theta_k, \quad \text{ for }k=0,\dots,n-1 .
		\]
		These equations imply the recursive relation $
		\lambda_k=\frac{1}{k+1} \Theta_{k-1} = \frac{-1}{k(k+1)}\lambda_{k-2},
		$ whose space of solutions is spanned by $\Psi_{\mathcal{E}_2,0}, \Psi_{\mathcal{E}_2,1}$. On the other hand, one easily checks that $\xi_2\cdot \Psi_{\mathcal{E}_2,0} = \Psi_{\mathcal{E}_2,1}$, hence $\Psi_{\mathcal{E}_2,0}, \Psi_{\mathcal{E}_2,1}$ are a basis for $\mathcal{E}_2$ as desired.
		\item \textul{$i=3$:} Similarly to the previous two subcases, we first note that
		\[
		\xi_3\cdot \psi_k = y_1\wedge \omega^{k+1} +i(k+1) \omega^k .
		\]
		Substituting this into (\ref{inducingKSMasterEqn}) gives
		\[
		\sum_{k=-1}^{n-1} \lambda_k[y_1\wedge \omega^{k+1} +i(k+1)\omega^k] =   \sum_{k=-1}^{n-1} \Theta_k[\omega^{k+1}-i(k+1)y_1\wedge \omega^k],
		\]
		which is equivalent to the linear system
		\[
		i(k+1)\lambda_k=\Theta_{k-1},\quad \lambda_{k-1}=-i(k+1)\Theta_k, \quad \text{ for }k=0,\dots,n-1 .
		\]
		These imply the recursive relation $
		\lambda_k = \frac{-i}{k+1} \Theta_{k-1} = \frac{1}{k(k+1)}\lambda_{k-2},
		$ whose space of solutions is spanned by $\Psi_{\mathcal{E}_3,0}$, $\Psi_{\mathcal{E}_3,1}$. The result in this subcase then follows by noting that $\xi_3\cdot \Psi_{\mathcal{E}_3,0} = i\Psi_{\mathcal{E}_3,1}$. 
	\end{enumerate}
	Thus, for $(M^{4n-1}, g)\ncong (S^{4n-1},g_{\text{round}})$ we have shown that $\mathcal{E}_i=E_i^-$ for $i=1,2,3$, and that $\{ \Psi_{\mathcal{E}_i,0},\Psi_{\mathcal{E}_i,1} \}$ is a basis (over $C^{\infty}(M)$) for this vector bundle. In particular, these basis spinors satisfy the algebraic conditions defining the bundles $E_i^-$. But the defining conditions for the $E_i^-$ are the same on $(M^{4n-1}, g)\cong (S^{4n-1}=\frac{\Sp(n)}{\Sp(n-1)},g_{\text{round}})$, and the above calculation of the basis spinors for $\mathcal{E}_i$ is still valid, hence the result also holds in this case. Supposing now that $n=1$, we have $E_i^- = \Sigma M$ for all $i=1,2,3$, and the result follows by noting that $\psi_{-1},\psi_0\in \cap_{i=1}^3 \mathcal{E}_i$ due to the relations
	\[\xi_1 \cdot \psi_{-1} = i \psi_{-1}, \qquad \xi_1\cdot \psi_0 = -i\psi_0,\qquad \xi_2\cdot \psi_{-1} = i\psi_0,\qquad \xi_3\cdot \psi_{-1} = \psi_0. \]    
\end{proof}

\begin{remark}
The basis $\Psi_{\mathcal{E}_1,0} = 1$, $\Psi_{\mathcal{E}_1,1} = y_1\wedge \omega^{n-1}$ for $E_1^- = \mathcal{E}_1$ was already known in \cite[Prop.\@ 4.15]{AHLspheres}, where it was computed directly using the defining condition of $E_1^-$; this calculation takes advantage of the relatively simple behaviour of $\omega^k$ and $y_1\wedge \omega^k$ under Clifford multiplication by $\xi_1$ and $e_i\cdot \varphi_1(e_i)$ (where $\{e_i\}_{i=1}^{4n-1}$ is an adapted basis). The corresponding Clifford products for the other two Sasakian structures are more complicated, and, consequently, direct calculations for $E_2^-$ and $E_3^-$ would require solving a difficult linear system where the number of equations grows with the dimension of the manifold. The preceding proposition avoids this problem by re-characterizing the bundles $E_i^-$ in terms of a desirable geometric property and solving the much easier equations resulting from this description.
\end{remark}

As noted in the discussion preceding Proposition \ref{Chap_Hof22:E_spaces_bases_hom_ts}, these spinors are precisely the ones which recover the homogeneous $\ts$ structure via the construction in Section \ref{constructionsection}, giving a full picture of this construction in the homogeneous $\ts$ setting:
\begin{theorem}
	If $(M=G/H,g,\xi_i,\eta_i,\varphi_i)$ is a simply-connected homogeneous $\ts$ space then, for each $i\in \{1,2,3\}$, the Sasakian structure $(\xi_i,\eta_i,\varphi_i)$ arises from the Killing spinors $\Psi_i:= \Psi_{\mathcal{E}_i,0}$ and $\Psi_i':= -\xi_i\cdot \Psi_{\mathcal{E}_i,0}$ (or $\Psi_i:= \Psi_{\mathcal{E}_i,1}$ and $\Psi_i':= -\xi_i\cdot \Psi_{\mathcal{E}_i,1}$) via the construction in Section \ref{constructionsection}.
\end{theorem}

The values of the spinors $\Psi_{\mathcal{E}_i,0}$ and $\Psi_{\mathcal{E}_i,1}$, in terms of the basis of invariant spinors from Theorem \ref{invspinortheorem}, are tabulated for a few low dimensions in Table \ref{Tab:spinors_PsiEi_low_dim}.

\begin{table}[h!]
	\centering
	\begin{tabular}{ |l||l|l|l| }
		\hline
		$\dim(M)$  &  $\Psi_{\mathcal{E}_1,0}$ & $\Psi_{\mathcal{E}_2,0}$ & $\Psi_{\mathcal{E}_3,0}$ \\
		\hline \hline 
		$7$   & $1$       &  $\omega-iy_1$   &  $\omega-iy_1$ \\
		$11$  & $1$     & $\omega-iy_1+\frac{1}{2}i y_1\wedge \omega^2$  & $\omega -iy_1 - \frac{1}{2}iy_1\wedge \omega^2$ \\
		$15$  & $1$      & $\omega-iy_1 +\frac{1}{2}i y_1\wedge \omega^2 -\frac{1}{6}\omega^3$ & $\omega -iy_1 - \frac{1}{2}iy_1\wedge \omega^2 +\frac{1}{6}\omega^3$ \\
		\hline\hline 
		$\dim(M)$ & $\Psi_{\mathcal{E}_1,1}$ & $\Psi_{\mathcal{E}_2,1}$ & $\Psi_{\mathcal{E}_3,1}$ \\ \hline \hline 
		$7$   & $y_1\wedge \omega$       & $1+iy_1\wedge \omega$ & $1-iy_1\wedge \omega$ \\
		$11$  & $y_1\wedge \omega^2$      & $1 + iy_1\wedge \omega   -\frac{1}{2}\omega^2 $  &  $1 - iy_1\wedge \omega +\frac{1}{2}\omega^2 $ \\
		$15$  & $y_1\wedge \omega^3$      &  $1 +iy_1\wedge \omega -\frac{1}{2}\omega^2  - \frac{1}{6} iy_1\wedge \omega^3$ & $1 -iy_1\wedge \omega +\frac{1}{2}\omega^2  - \frac{1}{6} iy_1\wedge \omega^3$\\
		\hline
	\end{tabular}
	\caption{The Spinors $\Psi_{\mathcal{E}_i,0}$ and $\Psi_{\mathcal{E}_i,1}$ in Low Dimensions}
	\label{Tab:spinors_PsiEi_low_dim}
\end{table}

\begin{remark}
	Table \ref{Tab:spinors_PsiEi_low_dim} emphasizes a subtle point about the relationship between individual Killing spinors and Einstein-Sasakian structures: a given Killing spinor $\psi$ may be used in combination with different choices of second Killing spinor $\psi'$ to yield different Sasakian structures. For example, in dimension $7$ the Killing spinors $\psi:= \Psi_{\mathcal{E}_2,1} = 1+iy_1\wedge \omega $ and $\psi':= -\xi_1\cdot \psi = -i-y_1\wedge \omega$ induce the Sasakian structure $(\xi_1,\eta_1,\varphi_1)$ by the construction in Section \ref{constructionsection}. On the other hand, the Sasakian structure induced by $\psi$ and $\psi'':= -\xi_2\cdot \psi = -i y_1 + \omega $ is $(\xi_2,\eta_2,\varphi_2)$. This behaviour stems from the fact that $\psi \in \Gamma(E_1^-)\cap \Gamma(E_2^-)$, so the spinor $\psi$ encodes information related to both of these Sasakian structures. Furthermore these two Sasakian structures generate $(\xi_3,\eta_3,\varphi_3)$ (see e.g.\@ \cite[p.\@ 556]{Fried90}), so in dimension $7$ the data of a homogeneous $\ts$ structure can be recovered from just three Killing spinors instead of the four assumed in Theorem \ref{spinorto3structure}. This improvement can be extended to other dimensions by noting that the spinors $\Psi_1:= \Psi_{\mathcal{E}_1,0} =1$ and $\Psi_1':=  -\xi_1\cdot \Psi_{\mathcal{E}_1,0} = -i\Psi_1 $ are linearly dependent and recover $(\xi_1,\eta_1,\varphi_1)$, so the data encoded by the three Killing spinors $\Psi_1$, $\Psi_2:=\Psi_{\mathcal{E}_2,0}$, $\Psi_2':= -\xi_2\cdot \Psi_2$ is sufficient to recover the whole homogeneous $\ts$ structure. Indeed, this observation is a manifestation of the fact that a (connected, simply-connected) manifold of dimension $4n-1$ carrying at least three linearly independent Killing spinors is necessarily $\ts$ (the general result is obtained in \cite{Bar}, and the earlier work \cite{Fried90} treats the special case of dimension $7$). Moreover, this lower bound cannot be improved further since there exist Einstein-Sasakian manifolds (which admit two linearly independent Killing spinors by \cite[Thm.\@ 1]{Fried90}) that do not carry a $\ts$ structure (see e.g.\@ \cite[Table 2]{nearly_parallel_g2}). 
\end{remark}

%% file: Steering_TeX_File.bbl
\begin{thebibliography}{BFGK91}

\bibitem[ACFH15]{dim67}
Ilka Agricola, Simon~G. Chiossi, Thomas Friedrich, and Jos Höll.
\newblock Spinorial description of $\text{SU}(3)$- and $\text{G}_2$-manifolds.
\newblock {\em Journal of Geometry and Physics}, 98:535–555, Dec 2015.

\bibitem[AD20]{3str}
Ilka Agricola and Giulia Dileo.
\newblock Generalizations of 3-{S}asakian manifolds and skew torsion.
\newblock {\em Advances in Geometry}, 20(3):331--374, 2020.

\bibitem[ADS21]{hom3alphadelta}
Ilka Agricola, Giulia Dileo, and Leander Stecker.
\newblock Homogeneous non-degenerate 3-$(\alpha,\delta)$-{S}asaki manifolds and
  submersions over quaternionic {K}ähler spaces.
\newblock {\em Annals of Global Analysis and Geometry}, 60(1):111–141, Apr
  2021.

\bibitem[AF10]{3Sasdim7}
Ilka Agricola and Thomas Friedrich.
\newblock 3-{S}asakian manifolds in dimension seven, their spinors and
  ${G}_2$-structures.
\newblock {\em Journal of Geometry and Physics}, 60(2):326–332, Feb 2010.

\bibitem[AH23a]{AHduality}
Ilka Agricola and Jordan Hofmann.
\newblock $\mathcal{H}$-{Killing} spinors and spinorial duality for homogeneous
  3-$(\alpha,\delta)$-{Sasaki} manifolds, 2023.
\newblock Preprint. \url{https://arxiv.org/abs/2309.16610}.

\bibitem[AH23b]{AHprojective}
Diego {Artacho} and Jordan Hofmann.
\newblock Invariant spinors on homogeneous projective spaces (in preparation),
  2023.

\bibitem[AHL23]{AHLspheres}
Ilka Agricola, Jordan Hofmann, and Marie-Amélie Lawn.
\newblock Invariant spinors on homogeneous spheres.
\newblock {\em Differential Geometry and its Applications}, 89:102014, 2023.

\bibitem[ANT23]{ANT_book_principal_fibre_bundles}
Ilka Agricola, Henrik Naujoks, and Marvin Theiss.
\newblock {\em Geometry of Principal Fibre Bundles (to appear)}.
\newblock 2023.

\bibitem[Arv03]{Arvan}
Andreas Arvanitoyeorgos.
\newblock {\em An Introduction to {L}ie Groups and the Geometry of Homogeneous
  Spaces}.
\newblock American Mathematical Society, 2003.

\bibitem[B{\"a}r93]{Bar}
Christian B{\"a}r.
\newblock Real {K}illing spinors and holonomy.
\newblock {\em Communications in Mathematical Physics}, 154:509--521, 1993.

\bibitem[BFGK91]{BFGK}
Helga Baum, Thomas Friedrich, Ralf Grunewald, and Ines Kath.
\newblock {\em Twistors and {K}illing Spinors on {R}iemannian Manifolds}.
\newblock B. G. Teubner Verlagsgesellschaft, 1991.

\bibitem[BG99]{BG_3Sas_paper}
Charles Boyer and Krzysztof Galicki.
\newblock 3-{Sasakian} manifolds.
\newblock In {\em Surveys in differential geometry. Vol. VI: Essays on Einstein
  manifolds. Lectures on geometry and topology, sponsored by Lehigh
  University's Journal of Differential Geometry}, pages 123--184. Cambridge,
  MA: International Press, 1999.

\bibitem[BG08]{SasakianGeometry}
Charles~P. Boyer and Krzysztof Galicki.
\newblock {\em Sasakian Geometry}.
\newblock Oxford Mathematical Monographs. Oxford University Press, 2008.

\bibitem[BGM94]{BG3Sas}
Charles~P. Boyer, Krzysztof Galicki, and Benjamin~M. Mann.
\newblock The geometry and topology of 3-{S}asakian manifolds.
\newblock {\em Journal f\"ur die reine und angewandte Mathematik},
  455:183--220, 1994.

\bibitem[Bla76]{blair_contact_manifolds_book}
David~E. Blair.
\newblock {\em Contact manifolds in {Riemannian} geometry}, volume 509 of {\em
  Lect. Notes Math.}
\newblock Springer, Cham, 1976.

\bibitem[CG88]{Cahen_Gutt_spin_struct_symmetric_spaces}
Michel Cahen and Simone Gutt.
\newblock Spin structures on compact simply connected {Riemannian} symmetric
  spaces.
\newblock {\em Simon Stevin}, 62(3-4):209--242, 1988.

\bibitem[DKL22]{invariantspinstructures}
Jordi {Daura Serrano}, Michael Kohn, and Marie-Am{\'e}lie Lawn.
\newblock G-invariant spin structures on spheres.
\newblock {\em Ann. Global Anal. Geom.}, 62(2):437--455, 2022.

\bibitem[DOP20]{homdata}
Cristina Draper, Miguel Ortega, and Francisco~J. Palomo.
\newblock Affine connections on 3-{Sasakian} homogeneous manifolds.
\newblock {\em Math. Z.}, 294(1-2):817--868, 2020.

\bibitem[FH91]{FultonHarris1991}
William Fulton and Joe Harris.
\newblock {\em Representation Theory: A First Course}.
\newblock Graduate Texts in Mathematics. Springer New York, 1991.

\bibitem[FK88]{FK88_french}
Thomas Friedrich and Ines Kath.
\newblock Vari{\'e}t{\'e}s {r}iemanniennes compactes de dimension 7 admettant
  des spineurs de {K}illing.
\newblock {\em C. R. Acad. Sci., Paris, S{\'e}r. I}, 307(19):967--969, 1988.

\bibitem[FK89]{dim5_Killing_spinors}
Thomas Friedrich and Ines Kath.
\newblock {Einstein manifolds of dimension five with small first eigenvalue of
  the Dirac operator}.
\newblock {\em Journal of Differential Geometry}, 29(2):263 -- 279, 1989.

\bibitem[FK90]{Fried90}
Thomas Friedrich and Ines Kath.
\newblock 7-dimensional compact {R}iemannian manifolds with {K}illing spinors.
\newblock {\em Communications in Mathematical Physics}, 133:543--561, 1990.

\bibitem[FKMS97]{nearly_parallel_g2}
Thomas Friedrich, Ines Kath, Andrei Moroianu, and Uwe Semmelmann.
\newblock On nearly parallel {{\(G_2\)}}-structures.
\newblock {\em J. Geom. Phys.}, 23(3-4):259--286, 1997.

\bibitem[{Fri}80]{Friedrich80}
Thomas {Friedrich}.
\newblock {Der erste {E}igenwert des {D}irac-{O}perators einer kompakten,
  {R}iemannschen {M}annigfaltigkeit nichtnegativer {S}kalarkr\"ummung}.
\newblock {\em {Math. Nachr.}}, 97:117--146, 1980.

\bibitem[Fri98]{Friedrich_submanifold}
Thomas Friedrich.
\newblock On the spinor representation of surfaces in {E}uclidean 3-space.
\newblock {\em Journal of Geometry and Physics}, 28(1-2):143--157, Nov 1998.

\bibitem[{Fri}00]{FriedrichBook}
Thomas {Friedrich}.
\newblock {\em {D}irac operators in {R}iemannian geometry}, volume~25.
\newblock Providence, RI: American Mathematical Society (AMS), 2000.

\bibitem[GRS23]{Leander_roots}
Oliver Goertsches, Leon Roschig, and Leander Stecker.
\newblock Revisiting the classification of homogeneous 3-{Sasakian} and
  quaternionic {K{\"a}hler} manifolds.
\newblock {\em Eur. J. Math.}, 9(1):28, 2023.
\newblock Id/No 11.

\bibitem[GW09]{GoodmanWallach}
Roe Goodman and Nolan~R. Wallach.
\newblock {\em Symmetry, Representations, and Invariants}.
\newblock Springer, 2009.

\bibitem[HN12]{physics_KS_paper}
Derek Harland and Christoph N{\"o}lle.
\newblock Instantons and {Killing} spinors.
\newblock {\em J. High Energy Phys.}, 2012(3):38, 2012.
\newblock Id/No 082.

\bibitem[HS90]{Hirzebruch_Slodowy_elliptic_genera}
Friedrich Hirzebruch and Peter Slodowy.
\newblock Elliptic genera, involutions, and homogeneous spin manifolds.
\newblock {\em Geom. Dedicata}, 35(1-3):309--343, 1990.

\bibitem[Kas71]{Kashiwada}
Toyoko Kashiwada.
\newblock A note on a {Riemannian} space with {Sasakian} 3-structure.
\newblock {\em Nat. Sci. Rep. Ochanomizu Univ.}, 22:1--2, 1971.

\bibitem[Kat00]{kath_Tduals}
Ines Kath.
\newblock Pseudo-{Riemannian} {{\(T\)}}-duals of compact {Riemannian}
  homogeneous spaces.
\newblock {\em Transform. Groups}, 5(2):157--179, 2000.

\bibitem[Kir86]{KirchbergKlaus-Dieter1986Aeft}
Klaus-Dieter Kirchberg.
\newblock An estimation for the first eigenvalue of the {D}irac operator on
  closed {K}\"{a}hler manifolds of positive scalar curvature.
\newblock {\em Annals of global analysis and geometry}, 4(3):291--325, 1986.

\bibitem[KN63]{KN1}
Shoshichi Kobayashi and Katsumi Nomizu.
\newblock {\em Foundations of differential geometry. {I}}, volume~15 of {\em
  Intersci. Tracts Pure Appl. Math.}
\newblock Interscience Publishers, New York, NY, 1963.

\bibitem[KN69]{KN2}
Shoshichi Kobayashi and Katsumi Nomizu.
\newblock {\em Foundations of Differential Geometry, Volume II}.
\newblock Interscience Publishers, 1969.

\bibitem[Kon75]{konishi75}
Mariko Konishi.
\newblock {On manifolds with Sasakian $3$-structure over quaternion Kaehler
  manifolds}.
\newblock {\em Kodai Mathematical Seminar Reports}, 26(2-3):194 -- 200, 1975.

\bibitem[Kuo70]{3Sas_structure_reduction}
{Ying-yan} Kuo.
\newblock {On almost contact $3$-structure}.
\newblock {\em Tohoku Mathematical Journal}, 22(3):325 -- 332, 1970.

\bibitem[LCL88]{LiE}
{Marc A. A. van} Leeuwen, Arjeh~M. Cohen, and Bert Lisser.
\newblock Li{E}: a computer algebra package for {L}ie group computations.
\newblock \url{http://wwwmathlabo.univ-poitiers.fr/~maavl/LiE/}, 1988.

\bibitem[LM89]{LM}
Herbert~Blaine {Lawson, Jr.} and Marie{-}Louise Michelsohn.
\newblock {\em Spin Geometry}.
\newblock Princeton University Press, 1989.

\bibitem[MS43]{MontgomerySamelson43}
Deane Montgomery and Hans Samelson.
\newblock Transformation groups of spheres.
\newblock {\em Annals of Mathematics}, 44(3):454--470, 1943.

\bibitem[MS14a]{GKSEinstein}
Andrei Moroianu and Uwe Semmelmann.
\newblock Generalized {K}illing spinors on {E}instein manifolds.
\newblock {\em International Journal of Mathematics}, 25(4):1450033--19, 2014.

\bibitem[MS14b]{GKSspheres}
Andrei Moroianu and Uwe Semmelmann.
\newblock Generalized {K}illing spinors on spheres.
\newblock {\em Annals of global analysis and geometry}, 46(2):129--143, 2014.

\bibitem[Nom54]{Nomizumap}
Katsumi Nomizu.
\newblock Invariant affine connections on homogeneous spaces.
\newblock {\em American Journal of Mathematics}, 76(1):33--65, 1954.

\bibitem[Sch08]{tensorFFTs}
Alexander Schrijver.
\newblock Tensor subalgebras and first fundamental theorems in invariant
  theory.
\newblock {\em Journal of Algebra}, 319(3):1305--1319, 2008.

\bibitem[Sem03]{Sem03_conformal_killing_forms}
Uwe Semmelmann.
\newblock Conformal {Killing} forms on {Riemannian} manifolds.
\newblock {\em Math. Z.}, 245(3):503--527, 2003.

\bibitem[Udr69]{udriste69}
Constantin Udrişte.
\newblock Structures presque coquaternioniennes.
\newblock {\em Bulletin mathématique de la Société des Sciences
  Mathématiques de la République Socialiste de Roumanie}, 13
  (61)(4):487--507, 1969.

\bibitem[Wan58]{Wangconnections}
{Hsien{-}chung} Wang.
\newblock On invariant connections over a principal fibre bundle.
\newblock {\em Nagoya Mathematical Journal}, 13:1--19, 1958.

\bibitem[Wan89]{Wang}
Mckenzie Y.~K. Wang.
\newblock Parallel spinors and parallel forms.
\newblock {\em Annals of Global Analysis and Geometry}, 7:59--68, 01 1989.

\bibitem[Wol65]{wolf_spaces}
Joseph~A. Wolf.
\newblock Complex homogeneous contact manifolds and quaternionic symmetric
  spaces.
\newblock {\em Journal of Mathematics and Mechanics}, 14(6):1033--1047, 1965.

\end{thebibliography}
